\documentstyle{amsppt}
\newcount\mgnf\newcount\tipi\newcount\tipoformule\newcount\greco 
\tipi=2          
\tipoformule=0   

\global\newcount\numsec\global\newcount\numfor
\global\newcount\numapp\global\newcount\numcap
\global\newcount\numfig\global\newcount\numpag
\global\newcount\numnf
\global\newcount\numtheo

\def\SIA #1,#2,#3 {\senondefinito{#1#2}%
\expandafter\xdef\csname #1#2\endcsname{#3}\else
\write16{???? ma #1,#2 e' gia' stato definito !!!!} \fi}

\def \FU(#1)#2{\SIA fu,#1,#2 }

\def\etichetta(#1){(\veroparagrafo.\veraformula)%
\SIA e,#1,(\veroparagrafo.\veraformula) %
\global\advance\numfor by 1%
\write15{\string\FU (#1){\equ(#1)}}%
\write16{ EQ #1 ==> \equ(#1)  }}

\def\etichettat(#1){\veroparagrafo.\veratheorema:%
\SIA e,#1,{\veroparagrafo.\veratheorema} %
\global\advance\numtheo by 1%
\write15{\string\FU (#1){\thu(#1)}}%
\write16{ TH #1 ==> \thu(#1)  }}

\def\etichettaa(#1){(A\veraappendice.\veraformula)
 \SIA e,#1,(A\veraappendice.\veraformula)
 \global\advance\numfor by 1
 \write15{\string\FU (#1){\equ(#1)}}
 \write16{ EQ #1 ==> \equ(#1) }}
\def\getichetta(#1){Fig. \verafigura
 \SIA g,#1,{\verafigura}
 \global\advance\numfig by 1
 \write15{\string\FU (#1){\graf(#1)}}
 \write16{ Fig. #1 ==> \graf(#1) }}
\def\retichetta(#1){\numpag=\pgn\SIA r,#1,{\verapagina}
 \write15{\string\FU (#1){\rif(#1)}}
 \write16{\rif(#1) ha simbolo  #1  }}
\def\etichettan(#1){(n\verocapitolo.\veranformula)
 \SIA e,#1,(n\verocapitolo.\veranformula)
 \global\advance\numnf by 1
\write16{\equ(#1) <= #1  }}

\newdimen\gwidth
\gdef\profonditastruttura{\dp\strutbox}
\def\senondefinito#1{\expandafter\ifx\csname#1\endcsname\relax}
\def\BOZZA{
\def\alato(##1){
 {\vtop to \profonditastruttura{\baselineskip
 \profonditastruttura\vss
 \rlap{\kern-\hsize\kern-1.2truecm{$\scriptstyle##1$}}}}}
\def\galato(##1){ \gwidth=\hsize \divide\gwidth by 2
 {\vtop to \profonditastruttura{\baselineskip
 \profonditastruttura\vss
 \rlap{\kern-\gwidth\kern-1.2truecm{$\scriptstyle##1$}}}}}
\def\verapagina{
{\romannumeral\number\numcap}.\number\numsec.\number\numpag}}

\def\alato(#1){}
\def\galato(#1){}
\def\veroparagrafo{\number\numsec}\def\veraformula{\number\numfor}
\def\veraappendice{\number\numapp}
\def\verapagina{\number\pageno}\def\veranformula{\number\numnf}
\def\verafigura{{\romannumeral\number\numcap}.\number\numfig}
\def\verocapitolo{\number\numcap}\def\veranformula{\number\numnf}
\def\veratheorema{\number\numtheo}
\def\Eqn(#1){\eqno{\etichettan(#1)\alato(#1)}}
\def\eqn(#1){\etichettan(#1)\alato(#1)}
\def\TH(#1){{\etichettat(#1)\alato(#1)}}
\def\thv(#1){\senondefinito{fu#1}$\clubsuit$#1\else\csname fu#1\endcsname\fi} 
\def\thu(#1){\senondefinito{e#1}\thv(#1)\else\csname e#1\endcsname\fi}
\def\ver{\veroparagrafo}
\def\Eq(#1){\eqno{\etichetta(#1)\alato(#1)}}
\def\eq(#1){\etichetta(#1)\alato(#1)}
\def\Eqa(#1){\eqno{\etichettaa(#1)\alato(#1)}}
\def\eqa(#1){\etichettaa(#1)\alato(#1)}
\def\dgraf(#1){\getichetta(#1)\galato(#1)}
\def\drif(#1){\retichetta(#1)}

\def\eqv(#1){\senondefinito{fu#1}$\clubsuit$#1\else\csname fu#1\endcsname\fi}
\def\equ(#1){\senondefinito{e#1}\eqv(#1)\else\csname e#1\endcsname\fi}
\def\graf(#1){\senondefinito{g#1}\eqv(#1)\else\csname g#1\endcsname\fi}
\def\rif(#1){\senondefinito{r#1}\eqv(#1)\else\csname r#1\endcsname\fi}
\def\bib[#1]{[#1]\numpag=\pgn
\write13{\string[#1],\verapagina}}

\def\include#1{
\openin13=#1.aux \ifeof13 \relax \else
\input #1.aux \closein13 \fi}

\openin14=\jobname.aux \ifeof14 \relax \else
\input \jobname.aux \closein14 \fi
\openout15=\jobname.aux
\openout13=\jobname.bib


\ifnum\tipoformule=1\let\Eq=\eqno\def\eq{}\let\Eqa=\eqno\def\eqa{}
\def\equ{}\fi


{\count255=\time\divide\count255 by 60 \xdef\hourmin{\number\count255}
        \multiply\count255 by-60\advance\count255 by\time
   \xdef\hourmin{\hourmin:\ifnum\count255<10 0\fi\the\count255}}

\def\oramin{\hourmin }

\def\data{\number\day/\ifcase\month\or january \or february \or march \or
april \or may \or june \or july \or august \or september
\or october \or november \or december \fi/\number\year;\ \oramin}

\setbox200\hbox{$\scriptscriptstyle \data $}

\newcount\pgn \pgn=1
\def\foglio{\number\numsec:\number\pgn
\global\advance\pgn by 1}
\def\foglioa{A\number\numsec:\number\pgn
\global\advance\pgn by 1}

\footline={\rlap{\hbox{\copy200}}\hss\tenrm\folio\hss}


\global\newcount\numpunt

\magnification=\magstephalf
\baselineskip=16pt
\parskip=8pt

\voffset=2.5truepc
\hoffset=0.5truepc
\hsize=6.1truein
\vsize=8.4truein 
{\headline={\ifodd\pageno\rightheadline \else \leftheadline \fi}}
\def\rightheadline{\it  {tralala}\hfil\tenrm\folio}
\def\leftheadline{\tenrm \folio \hfil\it  {Section $\ver$}}

\def\a{\alpha}

\def\d{\delta}
\def\e{\epsilon}

\def\f{\phi}
\def\g{\gamma}

\def\l{\lambda}

\def\s{\sigma}
\def\t{\tau}

\def\o{\omega}
\def\D{\Delta}
\def\L{\Lambda}
\def\G{\Gamma}
\def\O{\Omega}
\def\S{\Sigma}

\def\del #1{\frac{\partial^{#1}}{\partial\l^{#1}}}

\def\1{{1\kern-.25em\roman{I}}}
\def\eu{{1\kern-.25em\roman{I}}}
\def\f1{{1\kern-.25em\roman{I}}}

\def\R{{\Bbb R}}  
\def\N{{\Bbb N}}  
\def\P{{\Bbb P}}  
\def\Z{{\Bbb Z}}  
\def\Q{{\Bbb Q}}  
\def\C{{\Bbb C}}  
\def\E{{\Bbb E}}  

\def\la{\langle}
\def\ra{\rangle} 

\def\del{\partial}


\let\cal=\Cal

\def\DD{{\cal D}}
\def\EE{{\cal E}}

\def\GG{{\cal G}}

\def\KK{{\cal K}}
\def\LL{{\cal L}}
\def\MM{{\cal M}}
\def\NN{{\cal N}}
\def\OO{{\cal O}}

\def\QQ{{\cal Q}}
\def\RR{{\cal R}}
\def\SS{{\cal S}}
\def\TT{{\cal T}}

\def\chap #1#2{\line{\ch #1\hfill}\numsec=#2\numfor=1\numtheo=1}

\def\ba{{\backslash}}
\def\sb{{\subset}}

\def\em{{\emptyset}}


\def\note#1{\footnote{#1}}

\def\frac#1#2{{#1\over #2}}

\def\text#1{\quad{\hbox{#1}}\quad}
\def\newpage{\vfill\eject}
\def\proposition #1{\noindent{\thbf Proposition #1}}

\def\theo #1{\noindent{\thbf Theorem {#1} }}

\def\lemma #1{\noindent{\thbf Lemma {#1} }}
\def\definition #1{\noindent{\thbf Definition {#1} }}

\def\corollary #1{\noindent{\thbf Corollary #1 }}
\def\proof{{\noindent\pr Proof: }}
\def\proofof #1{{\noindent\pr Proof of #1: }}
\def\endproof{$\diamondsuit$}
\def\remark{\noindent{\bf Remark: }}
\def\thanks{\noindent{\bf Acknowledgements: }}

\font\pr=cmbxsl10

\font\thbf=cmbxsl10 scaled\magstephalf

\font\ch=cmbx12

\font\it=cmti10
\font\bf=cmbx10


\overfullrule=0pt
\def\op{\operatorname}

\font\tit=cmbx12
\font\aut=cmbx12

\def\s{\char'31}
\centerline{\tit METASTABILITY AND LOW LYING SPECTRA }
\vskip.2truecm
\centerline{\tit IN REVERSIBLE MARKOV CHAINS}
\vskip.2truecm 
\vskip1.5truecm

\centerline{\aut Anton Bovier 
\note{ Weierstrass-Institut f\"ur Angewandte Analysis und Stochastik,
Mohrenstrasse 39, D-10117 Berlin,\hfill\break Germany.
 e-mail: bovier\@wias-berlin.de},
Michael Eckhoff\note{Institut f\"ur Mathematik, Universit\"at Potsdam,
 Am Neuen Palais 10, D-14469 Potsdam, Germany.\hfill\break
e-mail: meckhoff\@math.uni-potsdam.de},}
\centerline{\aut V\'eronique Gayrard\note{DMA, EPFL, CH-1021 Lausanne, 
Switzerland, and
Centre de Physique Th\'eorique, CNRS,
Luminy, Case 907, F-13288 Marseille, Cedex 9, France. email: Veronique.Gayrard\@epfl.ch},
Markus Klein\note {Institut f\"ur Mathematik, Universit\"at Potsdam,
Am Neuen Palais 10, D-14469 Potsdam, Germany.\hfill\break
e-mail: mklein\@.math.uni-potsdam.de}}
\vskip1cm
\vskip0.5truecm\rm
\def\s{\sigma}
\parskip=0pt
\noindent {\bf Abstract:} We study a  large class of reversible 
 Markov chains with 
discrete state space and transition matrix $P_N$. 
We define the notion of a set of {\it metastable points} 
as a subset of the state space $\G_N$ such that (i) this set is reached 
from any point $x\in \G_N$ without return to $x$ with probability
at least $b_N$, while (ii) for any two point $x,y$ in the metastable set, the 
probability  $T^{-1}_{x,y}$ to reach $y$ from $x$ without return to $x$ 
is smaller than $a_N^{-1}\ll b_N$. Under some  additional 
  non-degeneracy assumption, we show that in such a situation:
\item{(i)} To each metastable point corresponds a metastable state,
whose mean exit time can be computed precisely.
\item{(ii)} To each metastable point corresponds one simple eigenvalue
of $1-P_N$ which is essentially equal to 
the inverse mean exit time from this state. 
Moreover, these results imply very sharp uniform control of the deviation
of the probability distribution of metastable exit times from the exponential
distribution.
\parskip=8pt


 \vfill
$ {} $
\newpage
{\headline={\ifodd\pageno\rightheadline \else \leftheadline \fi}}
\def\rightheadline{\it  {Metastability and spectra}\hfil\tenrm\folio}
\def\leftheadline{\tenrm \folio \hfil\it  {Section $\ver$}}

\def\op{\operatorname}
\def\vep{\varepsilon}

\chap{1. Introduction}1

In a recent paper [BEGK] we have presented rather sharp estimates on 
metastable transition times, both on the level of their mean values,
their Laplace transforms, and their distribution, for a class of reversible 
Markov chains that may best be characterized as random walks in
multi-well potentials, and that arise naturally in the context of
Glauber dynamics for certain mean field models. These results allow
for a very precise control of the behaviour of such processes over
very long times. 

In the present paper we continue our investigation of metastability 
in Markov chains focusing however on the connection between {\it 
metastability and
spectral theory} while working in a more general abstract context. 
Relating metastability to spectral characteristics of the Markov generator
or transition matrix is in fact a rather old topic. First mathematical 
results go back at least as far as Wentzell [W] and Freidlin and
Wentzell [FW]. Freidlin and Wentzell relate the eigenvalues 
of the transition matrix of Markov processes with exponentially\hfill\break
 small
transition probabilities to exit times from ``cycles''; Wentzell 
has a similar result for the spectral gap in the case of certain diffusion 
processes. All these relations are on the level of logarithmic equivalence,
i.e. of the form $\lim_{\e\downarrow 0} \e \ln (\l^\e_i T^\e_i) =0$
where $\e$ is the small parameter, and $\l^\e_i, T^\e_i$ are
the eigenvalues, resp. exit times. For more recent results
of this type, see [M,Sc]. Rather recently, Gaveau and Schulman 
[GS] (see also [BK] for an interesting discussion)
 have developed a more general program to give a spectral {\it definition}
of metastability in a rather general setting of Markov chains with 
discrete state space. In their approach low lying
 eigenvalues are related to metastable time scales and the corresponding 
eigenfunctions are related to metastable states. This interesting approach
still suffers, however, from rather imprecise relations 
  between eigenvalues and time-scales, and eigenfunctions and 
states. 

In this paper we will put these notions on a mathematically clean 
and precise basis for a wide class of Markov chains $X_t$ with countable 
state space $\G_N$\note{We expect that this approach can be extended with 
suitable
modifications to processes with continuous state space. Work on this problem 
is in progress.}, indexed by some large parameter $N$. 
Our starting point will be the definition of a 
{\it metastable set} of points 
each of which is supposed to be a representative
of one {\it metastable state}, on a chosen time scale.
It is important that our approach allows to consider the case where the 
cardinality of  
$\MM_N$ depends on $N$. 
 The key idea behind our definition will be that it ensures that 
 the time it takes to visit the representative point once the process
enters a ``metastable state'' is very short compared to the lifetime
of the metastable state. Thus, observing the visits of the process at the 
metastable set suffices largely  to trace the history of the process. 
We will then show that (under certain conditions ensuring the simplicity of 
the low-lying spectrum)  the expected times of transitions from 
each such metastable point to ``more stable'' ones (this notion will be defined
precisely later) are {\it precisely} equal to the inverse of one eigenvalue 
(i.e. $T_i=\l_i^{-1} (1+o(1))$) 
and that the corresponding eigenfunction is 
essentially the indicator function of the {\it attractor } of the 
corresponding  metastable point. This relation between times and eigenvalues
can be considered as the analogue of a quantum mechanical 
``uncertainty principle''.
Moreover, we will give precise formulas
expressing these metastable transition times in terms of escape probabilities
and the invariant measure. Finally, we will  derive uniform convergence 
results for the probability distribution of these times to the exponential 
distribution. Let us note that one main clue to the precise uncertainty 
principle is that we consider {\it transition times} between metastable points,
rather than {\it exit times from domains}. In the existing literature,
the problem of transitions between states involving
the passage through some ``saddle point'' (or ``bottle neck'') 
is almost persistently avoided 
 (for reasons that we have pointed out in the introduction 
of [BEGK]), except in one-dimensional situations where special methods 
can be used (as mentioned e.g. in  the very recent paper [GM]). 
But  the passage through the saddle point 
has a significant impact on the transition time which in general can be 
neglected only on the level of logarithmic equivalence\note{E.g. the 
lack of precision in the relation $T_M=\OO(1/(1-(1-\l)^t))$ in [GS] is 
partly due 
to this fact.}. Our results here,
together with those in [BEGK], appear to be the first that systematically 
control these effects.        

Let us now introduce our setting.
We consider a discrete time\note{There is  no 
difficulty in applying our results to continuous time chains by using 
suitable embeddings.} and specify our Markov chains by their transition 
matrix $P_N$ whose elements $p_N(x,y)$, $x,y\in\G_N$
denote the one-step transition probabilities of the chain.
In this paper we focus on the  case where the chain is {\it reversible}  
\note{The case of irreversible Markov chains will
be studied in a forthcoming publication [EK].} 
with respect to some probability measure $\Q_N$ on $\G_N$.  
We will always be  interested in the case where the cardinality of $\G_N$ is 
finite but tends to infinity as $N\uparrow\infty$.
Intuitively, metastability corresponds to a situation
 where the state space $\G_N$ can be 
decomposed into a 
number of disjoint components each containing a state such that the
time to reach one of these states from anywhere is much smaller than the 
time it takes to travel between any two of these states. We will now 
make this notion
precise. Recall from [BEGK] 
the notation $\t^x_I$
for the first instance the chain starting in $x$ at time $0$ reaches the set
$I\subset\G_N$,
$$
\t^x_I\equiv \inf\left\{t>0: X_t\in I\big |X_0=x\right\}
\Eq(0.1)
$$
\newpage
\definition {\TH(D.1)} {\it A set $\MM_N\subset\G_N$ 
will be called a set of {\rm metastable
points}  
if it satisfies the following assumptions. 
For  finite positive 
constants $a_N$, $b_N$
such that,
for some sequence $\vep_N\downarrow 0$,
$a_N^{-1}\leq \vep_N b_N$  
it holds that
\item{(i)} For all $z\in\G_N$, 
$$
\P\left[\t^z_{\MM_N}\leq\t^z_z\right] \geq b_N
\Eq(0.2)
$$
\item{(ii)} For any $x\neq y\in\MM_N$, 
$$
\P\left[\t^x_y<\t^x_x\right]\leq a^{-1}_N
\Eq(0.3)
$$
}

We associate with each $x\in\MM_N$ its {\it local valley}
$$
A(x)\equiv \left\{z\in\G_N: \P\left[\t^z_x=\t^z_{\MM_N}\right]
=   \sup_{y\in\MM_N     } \P\left[\t^z_y=\t^z_{\MM_N}\right]\right\}
\Eq(0.5)
$$
We will set 
$$
R_x\equiv  \frac{\Q_N(x)}{\Q_N(A(x))}
\Eq(0.7)
$$
and 
$$
\eqalign{
r_N\equiv &\max_{x\in\MM_N} R_x\leq 1
\cr
c^{-1}_N\equiv&\min_{x\in\MM_N} R_x>0
}
\Eq(0.7bis)
$$
Note that the sets $A(x)$ are not necessarily disjoint. We will
however show later that the  set of points that belong
to more than one local valley has very small mass under $\Q_N$. 
The above conditions do not fix $\MM_N$ uniquely.  It will be reasonable
to choose $\MM_N$ always such that 
for all $x\in\MM_N$, 
$$
\Q_N(x) =\sup_{z\in A(x)}\Q_N(z)
\Eq(0.6)
$$

The quantities $\P\left[\t^x_I\leq\t^x_x\right]$, $I\subset \MM_N$
furnish crucial characteristics of the chain. We will therefore introduce some
special notation for them: for   $I\subset \MM_N$ and $x\in\MM_N\ba I$,
 set
$$
T_{x,I}\equiv\left( \P[\t^x_I\leq\t^x_x]\right)^{-1}
\Eq(0.60)
$$
and 
$$
T_I\equiv \sup_{x\in\MM_N\ba I} T_{x,I}
\Eq(0.61)
$$
Note that these quantities depend on $N$, even though this is
suppressed in the notation. 

For simplicity we will consider in this paper only chains that satisfy
an additional assumption
of {\it non-degeneracy}:

\definition {\TH(D.2)}{\it We say that the family of  Markov chains  is {\rm
generic} on the level of the set $\MM_N$, if  there exists a sequence
$\e_N\downarrow 0$, such that
\item{(i)} For all pairs $x,y\in\MM_N$, and any set $I\subset\MM_N\ba\{x,y\}$
either
$T_{x,I}\leq \e_N T_{y,I}$ or
$ T_{y,I}\leq \e_N T_{x,I}$.
\item{(ii)} There exists $m_1\in\MM_N$, s.t. for all $x\in\MM_N\ba m_1$,
$\Q_N(x)\leq\e_N\Q_N(m_1)$. 
}

We can now state our main results. We do this in 
a slightly simplified form; more precise statements, containing explicit
estimates of the error terms,  will be formulated in the
later sections. 

\theo{\TH(A.1)}{\it Consider a discrete time Markov chain with state space
$\G_N$, transition matrix $P_N$, and  metastable set $\MM_N$ (as defined in 
Definition 1.1). Assume that
the chain is generic on the level $\MM_N$ in the sense of Definition 
2.1.  Assume further that $r_N \vep_N |\G_N||\MM_N|\downarrow 0$, 
and $r_Nc_N \e_N\downarrow 0 $, as $N\uparrow\infty$.
  For every $x\in\MM_N$ set $\MM_N(x)\equiv \{y\in\MM_N:
\Q_N(y)>\Q_N(x)\}$, define the {\rm metastable exit time}
$t_x\equiv \t^x_{\MM_N(x)}$. Then
\item{(i)} For any $x\in\MM_N$,
$$
\E\, t_x =R^{-1}_xT_{x,\MM_N(x)} (1+ o(1))
\Eq(T.1)
$$
\item{(ii)} For any $x\in\MM_N$, there exists an eigenvalue $\l_x$ 
of $1-P_N$ that satisfies 
$$
\l_x=\frac 1{\E\, t_x}\left(1+o(1)\right)
\Eq(T.2)
$$
Moreover, there exists a constant $c>0$ such that for all $N$
$$
\s(1-P_N) \ba \cup_{x\in \MM_N} \l_x \subset (c b_N|\G_N|^{-1},1]
\Eq(T.3)
$$
(here $\s(1-P_N)$ denotes the spectrum of $1-P_N$).  
\item{(iii)} If $\phi_x$ denotes the right-eigenvector of $P_N$ 
corresponding to the eigenvalue $\l_x$, normalized so that $\phi_x(x)=1$, then
$$
\phi_x(y)=\cases \P[\t^y_x<\t^y_{\MM_N(x)}](1+o(1)),&\text {if}  
\P[\t^y_x<\t^y_{\MM_N(x)}] \geq \e_N\cr
O(\e_N),&\text{otherwise}\endcases
\Eq(T.5)
$$ 
\item{(iv)} For any $x\in\MM_N$, for all $t>0$,
$$
\P[t_x >t\E\,t_x]= e^{-t(1+o(1))}(1+o(1))
\Eq(T.4)
$$
}

\remark We will see that $\P[\t^y_x<\t^y_{\MM_N(x)}]$ is extremely close to 
one for all $y\in A(x)$, with the possible exception of some points for which
$\Q_N(y)\ll\Q_N(x)$. Therefore, the corresponding (normalized) left 
eigenvectors 
$\psi_x(y)\equiv \frac{\Q_N(y)\phi_x(y)}{\sum_{z\in\G_N}\Q_N(y)\phi_x(y)}$ 
are to very good approximation  
equal to the invariant measure conditioned on the valley $A(x)$. 
As the invariant measure $\Q_N$ conditioned on $A(x)$ can be reasonably 
identified with a metastable state,
this 
establishes in a precise way the relation between eigenvectors and metastable
distributions. Brought to a point, our theorem then says that the left 
eigenfunctions of $1-P_N$ are the metastable states, the corresponding 
eigenvalues the mean lifetime of these states which can be computed 
in terms of exit probabilities via \eqv(T.1), and that the lifetime of a 
metastable state is exponentially distributed. 

\remark Theorem \thv(A.1) actually holds under slightly weaker 
hypothesis than those stated in Definition \thv(D.2). Namely,
as will become clear in the proof given in Section 5, the non-degeneracy of the
quantities $T_{x,I}$ is needed only for certain sets $I$.
On the other hand, if these weaker conditions fail, 
the theorem will no longer be true in this simple form. Namely. 
in a situation
where  certain subsets $\SS_i\subset\MM_N$ are such that for all
$x\in \SS_i$, $T_{x,I}$ (for certain relevant sets $I$, see Section 5)
differ only by constant factors, the eigenvalues and eigenfunctions
corresponding to this set will have to be computed specially through a
finite dimensional, non-trivial diagonalisation problem. While this can
in principle be done on the basis of the methods presented here,
we prefer to stay within the context of the more transparent generic 
situation for the purposes of this paper. Even more interesting situations 
crating genuinely new effect occur when degenerate subsets of states
whose  cardinality tends to infinity with $N$ are present. While these
fall beyond the scope of the present paper, the tools provided here and 
in [BEGK] can still of use, as is shown in [BBG].

Let us comment on the general motivation behind the formulation of Theorem 
\thv(A.1). The theorem allows, in a very general setting, to 
reduce all relevant quantities governing the metastable behaviour of a Markov 
chain to the computation of the key parameters,  $T_{x,y}$ and $R_x$,
$x,y\in \MM_N$. The first point to observe is that these quantities
are in many situations rather easy to control with good precision.
In fact, control of $R_x$ requires only knowledge of the invariant measure. 
Moreover, the ``escape probabilities'', $T_{x,y}^{-1}$, 
are related by a factor $\Q_N(x)$ 
to the {\it Newtonian capacity} of the point 
$y$ relative to $x$ and thus  
satisfy a {\it variational principle} that allows to express them in terms
of certain constraint  minima of the Dirichlet form of the Markov chain 
in question. In [BEGK] we have shown how this well-known fact (see e.g.
[Li], Section 6) can be used to give very sharp estimates on these 
quantities for the discrete diffusion processes studied there. Similar ideas
may be used in a wide variety of situations (for another example,
see [BBG]); we remind the reader that the same variational 
representation is at the basis of the ``electric network''
method [BS]. Let us mention that our general obsession with sharp 
results is motivated mainly by applications to {\it disordered} models
there the transition matrix $P_N$ is itself a random variable. Fluctuation
effects on the long-time behaviour 
provoked by the disorder can then only be analysed if sharp estimates 
on the relevant quantities are available. For examples see [BEGK, BBG].

In fact, in the setting of [BEGK], i.e. a random walk on $(\Z/N)^d\cap \L$
with reversible measure $\Q_N(x)= \exp(-NF_N(x))$, where
$F_N$ is ``close'' to some smooth function $F$ with 
finite number of local  minima satisfying some additional genericity
requirements, and the  natural choice for   $\MM_N$ being  
the set of local minima 
of $F_N$, the key quantities of Theorem \thv(A.1) were estimated
as
$$
b_N\geq c N^{-1/2}
\Eq(0.62)
$$
$$
\eqalign{
r_N&\leq c N^{-d/2}\cr
c_N&\leq CN^{d/2}
}\Eq(0.63)
$$
$$
T_{x,y} =e^{\OO(1)}N^{-(d-2)/2} e^{N[F_N(z^*(x,y))-F_N(x)]}
\Eq(0.64)
$$
where $z^*(x,y) $ is the position of the saddle point  between $x$ and $y$. 
Moreover, under the genericity assumption of [BEGK], 
$$
\e_N\leq e^{-N^\a}
\Eq(65)
$$
for some $\a>0$. The reader will check that Theorem \thv(A.1), together
with the precisions detailed in the later sections,  provides
very sharp estimates on the low-lying eigenvalues of $1-P_N$ and 
considerably sharpens the estimates on the
distribution function of the metastable transition times given in [BEGK].

Let us note that  Theorem \thv(A.1) allows to get results under
much milder regularity assumptions on the functions $F_N$ then were assumed in 
[BEGK]; in particular,
it is clear that one can deal with situations where an unbounded number of
``shallow'' local minima is present. Most of such minima can 
 simply be ignored in the definition of the metastable set $\MM_N$ which then will 
take into account only sufficiently deep minima. This is an important point 
in many applications, e.g. to  spin glass like models 
(but also molecular dynamics, as discussed below),
 where the number of local 
minima is expected to be very large (e.g. $\exp(a N)$), while the 
metastable behaviour is dominated by much fewer ``valleys''. 
For a discussion from a physics point of view, see e.g. [BK].

A second motivation for Theorem \thv(A.1) is given by recent work of
 Sch\"utte et al. [S,SFHD]. There, a numerical method for the analysis of
metastable conformational states of macromolecules is proposed 
that relies on the numerical investigation of the Gibbs distribution 
for the molecular equilibrium state via a Markovian molecular dynamics
(on a discretized state space). The key idea of the approach is to 
replace the time-consuming full simulation of the chain by a numerical
computation of the low-lying spectrum and the corresponding eigenfunctions, 
and to deduce from here results on the metastable states and their life times.
Our theorem allows to rigorously justify these deductions in a quantitative 
way in a setting that is sufficiently general to incorporate their situations.
  
The remainder of this article is organized as follows. In Section 2
we recall some basic notions, and more importantly, show that the knowledge
of $T_{x,y}$ for all $x,y\in \MM_N$ is enough to estimate more general 
transition probabilities. As a byproduct, we will show the existence of a
natural ``valley-structure'' on the state space, and the existence of
a natural (asymptotic) ultra-metric on the set $\MM_N$. In Section 3 
we show how to estimate mean transition times. The key result will be 
Theorem \thv(LL.5) which will imply the first assertion of Theorem \thv(A.1).
In Section 4 we begin our investigation of the relation between spectra
and transition times. The key result there is a characterization of
parts of the spectrum of $(1-P_N)$  in terms of the roots of 
some non-linear equation involving certain Laplace transforms of
transitions times, as well as a representation 
of the corresponding eigenvectors in terms of such Laplace transforms.
This together with some analysis of the properties of these 
Laplace transforms and an upper bound, using a Donsker-Varadhan [DV]
argument, will give sharp two-sided  estimates on the first 
eigenvalue
of general Dirichlet operators in terms of mean exit times. These estimates
will furnish a crucial input for Section 5 where we will prove that
the low-lying eigenvalues of $1-P_N$ are very close to the principal eigenvalues
of certain Dirichlet operators $(1-P_N)^{\S_j}$, with suitably 
constructed exclusion sets $\S_j$. This will prove the second assertion 
of Theorem \thv(A.1). In the course of the proof we will 
also provide rather precise estimates on the corresponding eigenfunction.
In the last Section we use the spectral information obtained before
to derive, using Laplace inversion formulas, very sharp estimates 
on the probability distributions of transition times. These will in particular
imply the last assertion of Theorem \thv(A.1).

\thanks We would like to thank Christof Sch\"utte and his collaborators
for explaining their approach to conformational dynamics 
and very motivating discussions.

\newpage
\chap{2. Some notation and elementary facts.}2

In this section we collect some useful notations and a number of more or less
simple facts that we will come back to repeatedly.

The most common notion we will use are the stopping times $\t^x_I$
defined in \eqv(0.1).
To avoid having to distinguish cases where $x\in I$, it will sometimes be
convenient to use the alternative quantities
$$
\s_I^x\equiv\min\{t\geq 0\,:
X_t\in I\,|\, X_0=x\}
\Eq(1.6a)$$
that take the value $0$ if $x\in I$.

Our analysis is largely based on the study of 
Laplace transforms of transition times. 
For $I\sb\G_N$ we denote by $(P_N)^I$ the Dirichlet 
operator 
$$
(P_N)^I\equiv 
\1_{I^c}P_N\,:\,\1_{I^c}\R^{\G_N}\to \1_{I^c}\R^{\G_N},
\qquad I^c\equiv \G_N\ba I
\Eq(1.3.1)
$$
Since our Markov chains are reversible with respect to the 
measure $\Q_N$, the matrix $(P_N)^I$ 
is a symmetric operator  on $\1_{I^c}\ell^2(\G_N,\Q_N)$ and thus
$$||(P_N)^I||=
\max \{|\l|\,|\,\l\in\s((P_N)^{I})\}
\Eq(1.3.2)
$$
where $||\cdot||$ denotes the operator norm induced by 
$\1_{I^c}\ell^2(\G_N,\Q_N)$. 
For a point $x\in\G_{N}$, subsets 
$I,J\sb\G_{N}$ and $u\in\C$, 
$\Re(u)<-\log||(P_N)^{I\cup J}||$, we define 
$$
G_{I,J}^x(u)\equiv
\E\bigl[e^{u\t_{I}^x}
\1_{\t_{I}^x\leq \t_{J}^x}\bigr]
=\sum_{t=1}^\infty 
e^{ut}\P[\t^x_I=t\leq \t^x_J]
\Eq(1.3.3)
$$
and
$$
K_{I,J}^x(u)\equiv
\E\bigl[e^{u\s_{I}^x}\1_{
\s_{I}^x\leq \s_{J}^x}\bigr]=
\left\{
\matrix
G_{I,J}^x(u) & \text{for} x\notin I\cup J,\\
1 & \text{for} x\in I,\\
0 & \text{for} x\in J\ba I
\endmatrix\right.
\Eq(1.3.4)
$$ 
The Perron-Frobenius theorem applied to the positive 
matrix $(P_N)^I$ implies that $G_{I,J}^x(u)$ and 
$K_{I,J}^x(u)$ converge locally uniformly on their 
domain of definition, more precisely 
$$
-\log||(P_{N})^{I}||=
\sup\{u\in\R\,|\,K_{I,I}^{x}(u)
\op{\,\,exists\,\,for\,\,all\,\,}x\notin I\}
\Eq(1.3.4a)$$

We now collect a number of useful standard results 
that follow trivially from the strong Markov property and/or 
reversibility, for easy reference.

From the strong Markov property one gets:

\lemma{\TH(P.3)} {\it
Fix $I,J,L\sb\G_{N}$. Then for all 
$\Re(u)<-\log||(P_N)^{I\cup J}||$ 
$$
G_{I,J}^x(u)=G_{I\ba L,J\cup L}^x(u)+
\sum_{y\in L}
G_{y,I\cup J\cup L}^x(u)K_{I,J}^{y}(u),
\qquad x\in\G_{N}
\Eq(1.3.6)
$$
}

In the following we will adopt the (slightly awkward) notation
$
P_{N}F^{x}\equiv
\sum_{z\in\G_{N}}P_{N}(x,z)F^z
$
The following are  useful specializations of 
the foregoing  result which we state without proof:

\corollary{\TH(P.5)} {\it
Fix $I,J\sb\G_{N}$. Then for $x\in\G_N$
$$
e^{u}P_{N}K_{I,J}^{x}(u)=G_{I,J}^x(u),
\qquad x\in\G_{N}
\Eq(1.5.1)
$$
and 
$$
(1-e^{u}P_{N})\del_u{K}_{I,J}^x(u)=G_{I,J}^x(u),
\qquad x\notin I\cup J
\Eq(1.5.2)
$$
where $\del_u$ denotes differentiation w.r.t. $u$.}

The following {\it renewal equation} will be used heavily: 

\corollary{\TH(P.6)} {\it
Let $I\sb\G_{N}$. Then for all $x\notin I\cup y$ and
$\Re(u)<-\log||(P_N)^{I\cup y}||$
$$
G_{y,I}^x(u)=
\frac{G_{y,I\cup x}^x(u)}{
1-G_{x,I\cup y}^x(u)}
\Eq(1.6.1)
$$
}

finally, from reversibility of the chain one has 

\lemma{\TH(P.7)} {\it
Fix $x,y\in\G_{N}$ and $I\sb \G_{N}$. Then 
$$
\Q_{N}(x)G_{y,I\cup x}^x=
\Q_{N}(y)G_{x,I\cup y}^y
\Eq(1.7.1)
$$
}

The next few  Lemmata imply the existence of a nested valley structure 
and that the knowledge of the quantities $T_{x,y}$ and the invariant measure 
are enough to control all transition probabilities with sufficient precision.
The main result is an approximate 
ultra-metric triangle
inequality. Let us define (the capacity of $x$ relative to $y$) 
$E(x,y) =\Q_N(x) T^{-1}_{x,y}$. 
We will  show that

\lemma{\TH(LL.6)}{\it Assume that 
$y,m\in\G_N$ and $J\sb\G_N\ba y\ba m$ such that
for $0<\d<\frac 12$,  $E(m,J) \leq\d E(m,y)$. Then 
$$
\frac{1-2\d}{1-\d} \leq \frac{
E(m,J)}{E(y,J)}\leq \frac1{1-\d}
\Eq(4.7.1)
$$
}

\proof We first prove the upper bound. We write 
$$
\P[\t^{m}_{J}<\t^{m}_{m}]=
\sum_{x\in J}\frac{\Q_{N}(x)}{\Q_{N}(m)}
\P[\t^{x}_{m}<\t^{x}_{J}]
\Eq(4.7.2)
$$
Now 
$$
\P[\t^{x}_{m}<\t^{x}_{J}]
=
\P[\t^{x}_{m}<\t^{x}_{J},\t^{x}_{y}<\t^{x}_{J}]+
\P[\t^{x}_{m}<\t^{x}_{J\cup y}]
{\P[\t^{m}_{J}<\t^{m}_{y\cup m}]\over
\P[\t^{m}_{J\cup y}<\t^{m}_{m}]}
\Eq(4.7.3)
$$
Now by assumption,
$$
{\P[\t^{m}_{J}<\t^{m}_{y\cup m}]\over
\P[\t^{m}_{J\cup y}<\t^{m}_{m}]}
\leq
{\P[\t^{m}_{J}<\t^{m}_{m}]\over
\P[\t^{m}_{y}<\t^{m}_{m}]}
\leq
\d
\Eq(4.7.4)
$$
Inserting \eqv(4.7.4) into \eqv(4.7.3) we arrive
at
$$
\P[\t^{x}_{m}<\t^{x}_{J}]
\leq
\P[\t^{x}_{y}<\t^{x}_{J},\t^{x}_{m}<\t^{x}_{J}]+
\d\P[\t^{x}_{m}<\t^{x}_{J\cup y}]
\leq
\P[\t^{x}_{y}<\t^{x}_{J}]+
\d\P[\t^{x}_{m}<\t^{x}_{J}]
\Eq(4.7.5)
$$
Inserting this inequality into \eqv(4.7.2) implies
$$
\P[\t^{m}_{J}<\t^{m}_{m}]\leq
(1-\d)^{-1}
{\Q_{N}(y)\over\Q_{N}(m)}
\P[\t^{y}_{J}<\t^{y}_{y}]
\Eq(4.7.6)
$$
We now turn to the lower bound. 
We first show that the assumption  implies
$$
\P[\t^{y}_{J}<\t^{y}_{m}]<
\d(1-\d)^{-1}
\Eq(4.7.7)
$$
Namely,
$$\eqalign{
\P[\t^{m}_{J}<\t^{m}_{m}]
&\geq
\P[\t^{m}_{y}<\t^{m}_{J}<\t^{m}_{m}]
=
\P[\t^{m}_{y}<\t^{m}_{J\cup m}]
\P[\t^{y}_{J}<\t^{y}_{m}]
}
\Eq(4.7.8)
$$
But 
$$
\eqalign{
\P[\t^{m}_{y}<\t^{m}_{J\cup m}] =
&\P[\t^{m}_{y}<\t^{m}_{m}] -\P[\t^m_J<\t^m_y<\t^m_m]
\cr&
\geq \P[\t^{m}_{y}<\t^{m}_{m}]-\P[\t^m_J<\t^m_m]\cr
& \geq\P[\t^{m}_{y}<\t^{m}_{m}](1-\d)
}
\Eq(4.7.8a)
$$
where the last inequality follows from the assumption. Thus 
$$
\P[\t^{m}_{J}<\t^{m}_{m}] \geq \P[\t^{m}_{y}<\t^{m}_{m}]
\P[\t^{y}_{J}<\t^{y}_{m}] (1-\d)
\Eq(4.7.9)
$$
Solving this inequality for $\P[\t^{y}_{J}<\t^{y}_{m}]$, the assumption
yields \eqv(4.7.7).

We continue as in the proof of the upper bound and
write for $x\in J$, using \eqv(4.7.7),
$$
\eqalign{
\P[\t^{x}_{y}<\t^{x}_{J}]
&=
\P[\t^{x}_{y}<\t^{x}_{J},\t^{x}_{m}<\t^{x}_{J}]
+
\P[\t^{x}_{y}<\t^{x}_{J\cup m}]
\P[\t^{y}_{J}<\t^{y}_{m}]
\cr&\leq
\P[\t^{x}_{m}<\t^{x}_{J}]
+
\P[\t^{x}_{y}<\t^{x}_{J}]
\d(1-\d)^{-1}
}\Eq(4.7.10)
$$
proving
$$
\P[\t^{x}_{y}<\t^{x}_{J}]\leq
\P[\t^{x}_{m}<\t^{x}_{J}]
\frac {1-\d} {1-2\d}
\Eq(4.7.11)
$$
Inserting \eqv(4.7.11) into \eqv(4.7.2) for $m\equiv y$
and, using once more \eqv(4.7.2) in the resulting
estimate, we obtain
$$
\P[\t^{y}_{J}<\t^{y}_{y}]\leq\frac
{1-\d}{ 1-2\d}
\frac{\Q_{N}(m)}{\Q_{N}(y)}
\P[\t^{m}_{J}<\t^{m}_{m}]
\Eq(4.7.12)
$$
which yields the lower bound in \eqv(4.7.1).
\endproof

\corollary {\TH(LL.88)} {\it Assume that $x,y,z\in\MM_N$. Then 
$$
E(x,y)\geq \frac 13\min\left(E(x,z),E(z,y)\right)
\Eq(U.1)
$$
}
\proof  By contradiction. Assume that
$E(x,y)< \frac 13 \min\left(E(x,z),E(z,y)\right)$. 
Then $E(x,y)<\frac 13E(x,z)$, and so by Lemma \thv(LL.6),
$$
\frac 12\leq \frac{E(x,y)}{E(z,y)}\leq \frac 32
\Eq(U.2)
$$
and in particular $E(y,z)\leq  2E(x,y)$, in contradiction with the 
assumption. \endproof

If we set 
$$
e(x,y)\equiv \cases -\ln E(x,y),&\text{if} x\neq y\cr
0,&\text {if} x=y
\endcases
\Eq(U.5)
$$
then Lemma \thv(LL.6) implies that $e$ furnishes an ``almost''  ultra-metric, 
i.e. it holds that $e(x,y)\leq \max(e(x,z),e(z,y)) +\ln 3$ 
which will turn out to be a useful tool later.
We mention that in the case of discrete diffusions in potentials,
the quantities $e(x,y)$ are essentially $N$ times the heights of the 
essential saddles between points $x$ and $y$.

The appearance of a natural ultra-metric structure on the set of metastable 
states under our  minimal assumptions is  interesting in itself.

A simple corollary of Lemma \thv(LL.6) 
 shows that the notion of elementary valleys, $A(m)$,
is reasonable in the sense that ``few'' points may belong to more than one 
valley.

\lemma {\TH(LL.3)}{\it Assume that $x,m\in \MM_N$ and $y\in \G_N$. Then 
$$
\P[\t^y_m<\t^y_y]\geq \e\text{and} \P[\t^y_x<\t^y_y]\geq \e
\Eq(L.4)
$$
implies that 
$$
\Q_N(y)\leq 2\e^{-1} \Q_N(m)\P[\t^m_x<\t^m_m]
\Eq(L.5)
$$
}
We leave the easy proof to the reader.

\newpage

\chap{3. Mean transition times}3

In this chapter we will prove various estimates of 
conditioned transition times 
$\E[\t_{I}^{x}|\t_{I}^x\leq\t_{J}^x]$, where 
$I\cup J\sb\MM_N$. The control obtained is crucial 
for the investigation of the low lying spectrum 
in Chapters 4 and 5. 
In the particular setting of the paper [BEGK], essentially the same types of 
estimates have been proven. Apart from re-proving these in the more abstract 
setting we consider here, we also present entirely different proofs
that avoid the inductive structure of the proofs given in 
[BEGK]. Instead, it uses heavily a representation formula for the 
Green's function (which first appeared in Section 3, Eq. (3.12) of [BEGK]
\note{More recently, the same formula was rederived by Gaveau and 
Moreau [GM] also for the non-reversible case.}).
While the new proofs are maybe less intuitive from a probabilistic 
point of view, they are considerably simpler.

\theo{\TH(G.1)} {\it
Fix a nonempty, irreducible, proper subset 
$\O\sb\G_N$. Let $(1-P_N)^{\O^c}$ denote the Dirichlet 
operator with zero boundary conditions at $\O^c$.
Then the Green's function defined as 
$G_N^{\O^c}(x,y)\equiv
((1-P_N)^{\O^c})^{-1}\1_{y}(x)$, $x,y\in\O$, is 
given by
$$
G_N^{\O^c}(x,y)={\Q_N(y)\over\Q_N(x)}
{\P[\s_x^y<\t_{\O^c}^y]
\over
\P[\t_{\O^c}^x<\t_x^x]}
\qquad(x,y\in\O)
\Eq(2.2.4)
$$
}

\proof This theorem follows essentially from the proof of
Eq. (3.12) of [BEGK]. 
Using e.g. the maximum principle, it follows 
that $(1-P_N)^{\O^c}$ is invertible. From 
\eqv(1.5.1) we obtain, using \eqv(1.3.4),
$$
(1-P_N)^{\O^c}K^y_{x,\O^c}(0)=
\1_{x}(y)G_{\O^c,x}^x(0)
\qquad(x,y\in\O)
\Eq(2.2.11)
$$
This function serves as a fundamental solution 
and we compute for $x,y\in\O$, using the symmetry of 
$(1-P_{N})^{\O^c}$, 
$$
\eqalign{
\Q_N(x)G_{\O^c,x}^x(0) G_N^{\O^c}(x,y)
=&
\la (1-P_N)^{\O^c}K_{x,\O^c}^{(\cdot)}(0),
G_N^{\O^c}(\cdot,y)\ra_{\Q_N}
\cr=&
\la K_{x,\O^c}^{(\cdot)}(0),
(1-P_N)^{\O^c}G_N^{\O^c}(\cdot,y)\ra_{\Q_N}
\cr=&
\Q_N(y)K_{x,\O^c}^y(0)
}
\Eq(2.2.12)
$$
This proves \eqv(2.2.4).\endproof

\remark
Observe that \eqv(2.2.4) still makes sense for 
$x\in\O$ and $y\in\partial\O$, where we define 
the boundary $\partial I$ of a set $I\sb\G_{N}$ 
to be 
$$
\partial I\equiv\{x\in I^c\,|\,\exists y\in I\,:
\,P_{N}(y,x)>0\}
\Eq(2.2.9)$$
For such $x$ and $y$ reversibility \eqv(1.7.1) and 
the renewal relation \eqv(1.6.1) for $u\equiv 0$ and 
$I\equiv\O^c$ imply 
$$
G_{N}^{\O^c}(x,y)=
\P[\t_{y}^x=\t_{\O^c}^x]
\qquad(x\in\O,y\in\partial\O)
\Eq(2.2.10)$$

Based on Theorem \thv(G.1) we can derive an alternative representation of 
a particular $h$-transform of the 
Green's function with $h(y)=\P[\t^y_I\leq \t^y_J]$ 
that will prove useful in the sequel.

\proposition {\TH(G.2)}{\it For every nontrivial partition  
$I\cup J=\O^c$ such that $I$ and $J$ are not empty 
and $I\ba J$ communicates with $\O$ we have 
$$
\P[\t_{I}^x\leq \t_{J}^x]^{-1}
G_N^{\O^{c}}(x,y)
\P[\t_{I}^y\leq \t_{J}^y]
=
{\P[\s_{y}^x<\t_{I}^x|\t_{I}^x\leq\t_{J}^x]
\over 
\P[\t_{\O^{c}}^{y}<\t_{y}^{y}]}
\D_{\O^{c}}(x,y),
\qquad x,y\in\O
\Eq(2.2.5)
$$
where 
$$
\D_{\O^{c}}(x,y)\equiv
{\P[\t_{\O^{c}}^y<\t_{y}^y]
\P[\s_{\O^{c}\cup y}^x<\t_{x}^x]
\over
\P[\t_{\O^{c}}^x<\t_{x}^x]
\P[\s_{\O^{c}\cup x}^y<\t_{y}^y]},
\qquad x,y\in\O
\Eq(2.2.6)
$$
Furthermore, 
$$
\frac 13\leq \D_{\O^c}(x,y)\leq 3
\Eq(2.2.7)
$$
}

\proof \eqv(2.2.5) is  a straightforward calculation that uses 
the renewal equation \eqv(1.6.1), reversibility, and the strong Markov
property. Indeed, by \eqv(2.2.4) the left-hand side of \eqv(2.2.5)
equals
$$
\frac{\Q_N(y)\P[\s^y_x<\t^y_{\O^c}]\P[\t^y_I\leq \t^y_J]}
{\Q_N(x)\P[\s^x_{\O^c}<\t^x_x]\P[\t^x_I\leq \t^x_J]}
\Eq(2.1)
$$
By the renewal equation, this equals
$$
\frac{\Q_N(y)\P[\s^y_x<\t^y_{\O^c\cup y}]\P[\t^y_I\leq \t^y_J]}
{\Q_N(x)\P[\t^y_{\O^c\cup x}<\t^y_y]\P[\s^x_{\O^c}<\t^x_x]\P[\t^x_I\leq 
\t^x_J]}
\Eq(2.2)
$$
which by reversibility turns into
$$
\eqalign{
&\frac{\P[\s^x_y<\t^x_{\O^c\cup x}]\P[\t^y_I\leq \t^y_J]}
{\P[\t^y_{\O^c\cup x}<\t^y_y]\P[\s^x_{\O^c}<\t^x_x]\P[\t^x_I\leq \t^x_J]}
\cr&=
\frac{\P[\s^x_y<\t^x_{\O^c}]\P[\s^y_{\O^c\cup x}<\t^y_y]\P[\t^y_x<\t^y_{\O^c}]
\P[\t^y_I\leq \t^y_J]}
{\P[\t^y_{\O^c\cup x}<\t^y_y]\P[\s^x_{\O^c}<\t^x_x]\P[\t^x_I\leq \t^x_J]}
\cr
&=\frac{\P[\s^x_y<\t^x_I|\t^x_I\leq \t^x_J]\P[\s^y_{\O^c\cup x}<\t^y_y]\P[\t^y_x<\t^y_{\O^c}]}
{\P[\t^y_{\O^c\cup x}<\t^y_y]\P[\s^x_{\O^c}<\t^x_x]}
}
\Eq(2.3)
$$
where the last identity uses that by the strong Markov property
$$
\eqalign{
&\P[\s^x_y<\t^x_I,\t^x_I\leq \t^x_J]=\P[\s^x_y<\t^x_I\leq\t^x_J]
=\P[\s^x_y<\t^x_{I\cup J}]\P[\t^y_I\leq \t^x_J]
}
\Eq(2.4)
$$
\eqv(2.3) immediately implies \eqv(2.2.5).

We now turn to the proof of the bound \eqv(2.2.7). Since 
$\D_{\O^c}(x,x)=1$ it is enough to consider the case where $x\neq y$. 
Moreover, since $\D_{\O^c}(x,y)=\frac1{\D_{\O^c}(y,x)}$, an upper bound 
$\D_{\O^c}(x,y)\leq 3$ will immediately imply the claimed lower bound. 

The basic input here is the observation that a path from $y$ to $\O^c$ either 
visits a point $x$ or it does not, yielding, together with the strong Markov 
property
$$
\eqalign{
\P[\t^y_{\O^c}<\t^y_y]=&\P[\t^y_{\O^c}<\t^y_{x\cup y}]
+\P[\t^y_x<\t^y_{\O^c}<\t^y_{ y}]
\cr
&=\P[\t^y_{\O^c}<\t^y_{x\cup y}]
+\P[\t^y_x<\t^y_{\O^c\cup y}]\P[\t^x_{\O^c}<\t^x_{y}]
}
\Eq(2.5)
$$
Using this identity for the first factor  in the numerator of \eqv(2.2.6),
we obtain that 
$
\D_{\O^c}(x,y)$ can be written as $\D_{\O^c}(x,y)=(I)+(II)$
where
$$
\eqalign{
(I)=&\frac{\P[\t^y_{\O^c}<\t^y_{x\cup y}]\P[\t^x_{\O^c\cup y}<\t^x_x]}
{\P[\t^x_{\O^c}<\t^x_{x}]\P[\t^y_{\O^c\cup x}<\t^y_y]}
=\frac{\P[\t^y_{\O^c}<\t^y_{x}]\P[\t^x_{\O^c\cup y}<\t^x_x]}
{\P[\t^x_{\O^c}<\t^x_{x}]}
}
\Eq(2.6)
$$
The renewal equation was used in the second equality. Decompose the
event in the second factor of the numerator and use \eqv(2.5) in the 
denominator. This yields
$$
\eqalign{
(I)=&\frac{\P[\t^y_{\O^c}<\t^y_{x}]\left(\P[\t^x_{\O^c}<\t^x_{x\cup y}]
+\P[\t^x_y<\t^x_{\O^c\cup x}]\right)}
{\P[\t^x_{\O^c}<\t^x_{x\cup y}]+
\P[\t^x_y<\t^x_{\O^c\cup x}]\P[\t^y_{\O^c}<\t^y_{x}]
}
\leq \P[\t^y_{\O^c}<\t^y_{x}] +1\leq 2
}
\Eq(2.7)
$$
For (II) we get 
$$
\eqalign{
(II)&=\frac{\P[\t^y_x<\t^y_{\O^c\cup y}]\P[\t^x_{\O^c}<\t^x_{x\cup y}]
\P[\t^x_{\O^c\cup y}<\t^x_x]}
{\P[\t^x_{\O^c\cup y}<\t^x_x]
\P[\t^x_{\O^c}<\t^x_{x}]\P[\t^y_{\O^c\cup x}<\t^y_y]}
=\frac{\P[\t^y_x<\t^y_{\O^c\cup y}]\P[\t^x_{\O^c}<\t^x_{x\cup y}]}
{\P[\t^x_{\O^c}<\t^x_{x}]\P[\t^y_{\O^c\cup x}<\t^y_y]}
\leq 1
}
\Eq(2.8)
$$
The bounds \eqv(2.2.7) are now obvious. \endproof

The representation \eqv(2.2.5) for the Green's function implies immediately
a corresponding representation for the (conditioned) expectation
of entrance times $\t^x_I$. To see this, recall from \eqv(1.5.2) 
for $u\equiv 0$ that 
$$
(1-P_{N})^{I\cup J}
\E\left[\s_{I}^{y}\1_{\{\s_{I}^y\leq\s_{J}^y\}}\right]=
\P[\t_{I}^y\leq\t_{J}^y],
\qquad y\notin I\cup J
\Eq(2.3.2)
$$ 
This yields immediately the

\corollary{\TH(G.4)} {\it
Let $I,J\sb\G_N$. Then for all $x\notin I\cup J$
$$
\eqalign{
\E[\t_{I}^x|\t_{I}^x\leq\t_{J}^x]
&=\sum_{y\in (I\cup J)^c} \P[\t_{I}^x\leq \t_{J}^x]^{-1}
G_N^{\O^{c}}(x,y)
\P[\t_{I}^y\leq \t_{J}^y]\cr
&=
\sum_{y\in (I\cup J)^c}
{\Q_{N}(y)\over\Q_{N}(x)}
{\P[\s_x^y<\t_{I\cup J}^y]
\over
\P[\t_{I\cup J}^x<\t_x^x]}
{\P[\t_{I}^y\leq\t_{J}^y]
\over
\P[\t_{I}^x\leq\t_{J}^x]}
}
\Eq(2.3.1)
$$
}
 
A first consequence of the representation given 
above is  

\corollary{\TH(G.5)} {\it
Fix $I\sb\MM_{N}$. Then for all $x\in\G_{N}$
$$
\E[\t_{I}^x|\t_{I}^x<\t_{\MM_N\ba I}^x]\leq 3b_N^{-1} |\G_N|
\Eq(2.5.1)
$$
In particular, 
$$
\E[\t_{\MM_{N}}^x]\leq 3b^{-1}_N|\G_N|
\Eq(2.5.2)
$$
}

\proof Using \eqv(2.2.6) in \eqv(2.5.1), we get that
$$
\E[\t_{I}^x|\t_{I}^x<\t_{\MM_N\ba I}^x]
=\sum_{y\in \G_N\ba \MM_N} 
\frac{\P[\s^x_y<\t^x_I|\t^x_I\leq \t^x_{\MM_N\ba I}]}
{\P[\t^y_{\MM_N}<\t^y_y]} \D_{\MM_N}(x,y)
\Eq(2.9)
$$
Using the lower bound \eqv(0.2) from Definition 1.1 together with 
the upper bound \eqv(2.2.7), we get
$$
\E\left[\t_{I}^x|\t_{I}^x<\t_{\MM_N\ba I}^x\right] \leq 3 b_N^{-1}
 \sum_{y\in \G_N\ba \MM_N} \P[\s^x_y<\t^x_I|\t^x_I\leq \t^x_{\MM_N\ba I}]
\Eq(2.10)
$$
from which the claimed estimate follows by bounding the conditional
probability by one\note{It is obvious that in cases when $|\G_N|=\infty$
this bound can in many cases be improved to yield a reasonable estimate.
Details will however depend upon assumptions on the global geometry.}.
The special case $I=\MM_N$ follows in the same way, 
with the more explicit bound
$$
\E\t^x_{\MM_N} \leq 3b_N^{-1} \sum_{y\in\G_N\ba \MM_N} \P[\s^x_y<\t^x_{\MM_N}]
\Eq(2.11)
$$
This concludes the proof of the corollary.\endproof

Theorem 2.2 allows to compute  very easily 
the mean times of metastable transitions.

\theo{\TH(LL.5)}{\it  Assume that $J\subset\MM_N$, $x\in\MM_N$, and
$x,J$ satisfy the condition
$$
T_{x,J}=T_J
\Eq(L.21)
$$
Then 
$$
\E \t^x_J=\frac{\Q_N(A(x))}{\Q_N(x)\P[\t^x_J<\t^x_x]}
\left(1+ \OO(1)  \left( \frac{R_x |\MM_N| |\G_N|}{b_Na_N}+\e_N
R_x c_N\right)\right)
\Eq(L.11)
$$
}

\proof 
Specializing Corollary \thv(G.4) to the case $J=\emptyset$, we get the 
representation
$$
\E\t^x_J=\frac 1{\Q_N(x)\P[\t^x_J<\t^x_x]}
\sum_{y\not\in J} \Q_N(y)\P[\s^y_x<\t^y_J]
\Eq(L.12)
$$
We will decompose the sum into three pieces corresponding to the two
sets
$$\eqalign{
\O_1\equiv A(x)&\cr
\O_2\equiv & \G_N\ba A(x)\ba J}
\Eq(L.13)
$$
The sum over $\O_1$ gives the main contribution; the trivial upper bound
$$
\sum_{y\in\O_1} \Q_N(y) \P[\s^y_x<\t^y_J] \leq \sum_{y\in\O_1} \Q_N(y)
\Eq(L.15)
$$
is complemented by a lower bound that uses (we ignore the trivial case 
$x=y$ where\hfill\break $\P[\s^x_x<\t^x_J]=1$) 
$$
\P[\t^y_x<\t^y_J] =1-\P[\t^y_J<\t^y_x]\geq
1-\frac{\P[\t^y_J<\t^y_y]}{\P[\t^y_x<\t^y_y]}
\Eq(L.16)
$$
By Lemma \thv(LL.6), if
$\P[\t^x_J<\t^x_x]\leq \frac 13 \P[\t^x_y<\t^x_x]$, then  
$$
\P[\t^y_J<\t^y_y] \leq \frac 32
\frac{\Q_N(x)}{\Q_N(y)} \P[\t^x_J<\t^x_x] 
\Eq(L.17)
$$
so that 
$$
\Q_N(y)\frac{\P[\t^y_J<\t^y_y]}{\P[\t^y_x<\t^y_y]} \leq 
\frac 32\Q_N(x) \frac {|\MM_N|}{b_Na_N}
\Eq(L.17a)
$$
On the other hand, if $\P[\t^x_J<\t^x_x]> \frac 13 \P[\t^x_y<\t^x_x]$, then
$$
\Q_N(y)\leq 3\Q_N(x)\frac {\P[\t^x_J<\t^x_x]}{\P[\t^y_x<\t^y_y]}
\leq  3\Q_N(x)  \frac {|\MM_N|}{b_Na_N}
\Eq(L.17b)
$$
Thus
$$
\eqalign{
\sum_{y\in\O_1} \Q_N(y) \P[\s^y_x<\t^y_J] &\geq 
 \sum_{y\in\O_1} \Q_N(y) - 3|A(x)| \Q_N(x)   \frac {|\MM_N|}{b_Na_N}
\cr
&=\Q_N(A(x))\left(1- 3|A(x)| R_x  \frac {|\MM_N|}{b_Na_N}
\right)
}
\Eq(L.18)
$$
We now consider the remaining contributions. This is bounded by
$$
\frac 1{\Q_N(x)\P[\t^x_J<\t^x_x]}\sum_{m\in \MM\ba x} L_m
\Eq(L.18a)
$$
where
$$
L_m\equiv \sum_{y\in A(m)\ba J}L_m(y)\equiv
 \sum_{y\in A(m)\ba J}  \Q_N(y)  \P[\s^y_x<\t^y_J]
\Eq(L.18b)
$$
Assume first that $y$ is such that 
\item{(CJ)} $\Q_N(y) \P[\t^y_J<\t^y_y]\sim \Q_N(m)\P[\t^m_J<\t^m_m]$ and 
\item{(Cx)} $\Q_N(y) \P[\t^y_x<\t^y_y]\sim \Q_N(m)\P[\t^m_x<\t^m_m]$ hold,

where we introduced the notation 
 $a\sim b\Leftrightarrow \frac 13\leq \frac ab\leq 3$. 
Then 
$$
L_m(y) \leq 9\Q_N(y) \frac{\P[\t^m_x<\t^m_m]}{\P[\t^m_J<\t^m_m]}
\Eq(L.19)
$$
There are two cases:
\item{(i)} If $E(m,J)\leq \frac 13E(m,x)$, then by Lemma \thv(LL.6),
$
\frac{\Q_N(m)\P[\t^m_J<\t^m_m]}{ \Q_N(x)\P[\t^x_J<\t^x_x]}\leq \frac 32	
$
or 
$$
\Q_N(m)\leq \frac 32\Q_N(x) \frac {T_{m,J}}{T_{x,J}}\leq \e_N \frac 32\Q_N(x)
\Eq(L.200)
$$
Hence 
$$
L_m(y) \leq \Q_N(y) \leq \frac {\Q_N(y)}{\Q_N(m)} \e_N\frac 32
R_x\Q_N(A(x))
\Eq(L.201)
$$
\item {(ii)}   If $E(m,J) >  \frac 13E(m,x)$, then 
$E(x,J)\geq \frac 13E(m,x)$ or
$
\Q_N(x)\P[\t^x_J<\t^x_x]\geq \frac 13\Q_N(m)\P[\t^m_x<\t^m_m]
$
so that
$$
L_m(y)\leq 27\frac{\Q_N(y) \Q_N(x)}{\Q_N(m)}\frac{T_{m,J}}{T_{x,J}}
\leq 27\e_N R_x\frac{\Q_N(y) }{\Q_N(m)} \Q_N(A(x))
\Eq(L.202)
$$
Finally we must consider the cases where (CJ) or (Cx) are violated.
\item{(iii)} Assume that (Cx) fails. Then by Lemma \thv(LL.6),
$\P[\t^m_x<\t^m_m]\geq \frac 13 \P[\t^m_y<\t^m_m]$ which implies that
$$
\eqalign{
L_m(y)&\leq \Q_N(y)\leq 3\Q_N(m) \frac {\P[\t^m_x<\t^m_m]}{\P[\t^y_m<\t^y_y]}
\leq 3\Q_N(m) \frac {\P[\t^m_x<\t^m_m] |\MM_N|}{b_N}
	\cr
& \leq3\Q_N(x) \frac {\P[\t^x_m<\t^x_x] |\MM_N|}{b_N}\leq 		
\frac {3|\MM_N|}{b_N a_N}   R_x \Q_N(A(x))
}
\Eq(L.203)
$$
  \item{(iv)} Finally it remains the case where (CJ) fails but (Cx) holds.
Then $\P[\t^y_J<\t^y_y]>\frac 13 \P[\t^y_m<\t^y_y]\geq \frac {b_N}{3|\MM_N|}$
and
$\Q_N(y)\P[\t^y_x<\t^y_y]\leq \frac 32 \Q_N(m)\P[\t^m_x<\t^m_m]
=\frac32\Q_N(x)\P[\t^x_m<\t^x_x]$. Thus $L_m(y)$ satisfies equally
the bound \eqv(L.203).

Using these four bounds, 
summing over $y$ one gets
$$
L_m\leq 27\Q_N(A(x))
\max\left(\e_NR_x R_m^{-1}, \frac{ |\MM_N||A(m)|}{b_N a_N}R_x \right)
\Eq(L.204)
$$
Putting everything together, we arrive at the assertion of the theorem.\endproof

\remark As a trivial corollary from the proof of Theorem \thv(LL.5) one has

\corollary {\TH(LL.6a)} {\it Let $x\in \MM_N$ and $J\subset\MM_N(x)$.
Then the conclusions of Theorem \thv(LL.5) also hold.}

Finally, we can easily prove a general upper bound on any conditional
expectation.

\theo{\TH(LL.7)} {\it For any $x\in\G_N$ and $I,J\subset\MM_N$,
$$
\E\left[\t^x_I|\t^x_I\leq \t^x_J\right]
\leq C  \sup_{m\in\MM_N\ba I\ba J} 
\left(R_m\P[\t^m_{I\cup J}<\t^m_m]\right)^{-1}
\Eq(L.35)
$$
}

To prove this theorem the representation of the Green's function
given in Proposition 2.2 is particularly convenient. It yields
$$
\eqalign{
\E\left[\t^x_I|\t^x_I\leq \t^x_J\right]&=\sum_{y\in\G_N\ba I\ba J}
\frac {\P[\s^x_y<\t^x_I|\t^x_I\leq \t^x_J]}{\P[\t^y_{I\cup J}<\t^y_y]}
\D_{I\cup J}(x,y)}
\Eq(L.36)
$$
Note first that the terms with $y$ such that $\P[\t^y_{I\cup J}<\t^y_y]
\geq \d b_N$ 
yield a contribution of no more than
$|\G_N|(\d b_N)^{-1}$ which is negligible. To treat the remaining 
terms, we use that whenever 
$y\in A(m)$, Lemma \thv(LL.6)  implies that
$\P[\t^y_{I\cup J}<\t^y_y]\geq \frac{\Q_N(m)}{\Q_N(y)} 
\P[\t^m_{I\cup J}<\t^m_m]
$.
Thus 
$$
\eqalign{
\E\left[\t^x_I|\t^x_I\leq \t^x_J\right]&\leq \frac{3|\G_N|}{\d b_N}
+\sum_{m\in\MM_N\ba I\ba J}\sum_{y\in A(m)} 3\frac{\Q_N(y)}{\Q_N(m)} 
\frac{\P[\s^x_y<\t^x_I|\t^x_I\leq \t^x_J]}{\P[\t^m_{J\cup I}<\t^m_m]}\cr
&\leq  \frac{3|\G_N|}{\d b_N}
+\sum_{m\in\MM_N\ba I\ba J} 3R^{-1}_m
\frac 1{\P[\t^m_{J\cup I}<\t^m_m]}
}
\Eq(L.37)
$$
from which the claim of the theorem follows by our general assumptions.
Note that by very much the same arguments as used before, it is possible to 
prove  that 
$$
\P[\s^x_y<\t^x_I|\t^x_I\leq \t^x_J]
\leq (1+\d) \P[\s^x_m<\t^x_I|\t^x_I\leq \t^x_J]
\Eq(L.38)
$$
which allows to get the sharper estimate
$$
\E\left[\t^x_I|\t^x_I\leq \t^x_J\right]\leq \frac{3|\G_N|}{\d b_N}
+\sum_{m\in\MM_N\ba I\ba J} 3(1+\d) R^{-1}_{m}
\frac{\P[\s^x_m<\t^x_I|\t^x_I\leq \t^x_J]}{\P[\t^m_{J\cup I}<\t^m_m]}
\Eq(L.39)
$$
\endproof

We conclude this section by stating some consequences of the two preceding 
theorems that will be useful later. 

\lemma{\TH(LL.8)} {\it
Let $I,m$ satisfy the hypothesis of Theorem \thv(LL.5). Then
$$
\max_{x\notin I}\E[\t_{I}^{x}]=
\E[\t_{I}^{m}]
\bigl(
1+\OO(T_{I\cup m}/T_I)
\bigr)
\Eq(2.6.0)
$$
Moreover, we have
$$
{\E[\t_{m}^{m},\t_{m}^{m}<\t_{I}^{m}]
\over
\P[\t_{I}^{m}<\t_{m}^{m}]}=
\E[\t_{I}^{m}]
\bigl(
1-\OO(T_{I\cup m}/T_I))
\bigr)
\Eq(2.6.1)
$$
In particular,
$$
\E[\t_{m}^{m},\t_{m}^{m}<\t_{I}^{m}]=
R^{-1}_m 
\left(1+\OO(T_{I\cup m}/T_I)\right)
\Eq(2.6.2)
$$
}

\proof Decomposing into the events where $m$ is and is not visited before
$I$, and,
using the strong Markov property, one gets  
$$
\E[\t_{I}^x]=
\P[\t_{I}^x<\t_{m}^x]
\E[\t_{I}^x|\t_{I}^x<\t_{m}^x]
+
\P[\t_{m}^x<\t_{I}^x]
\bigl(\E[\t_{m}^x|\t_{m}^x<\t_{I}^x]+
\E[\t_{I}^m]\bigr)
\Eq(2.6.0a)
$$
Using Theorems \thv(LL.5) and \thv(LL.7), this implies 
\eqv(2.6.0) readily. 
In the same way, or by differentiating the renewal equation \eqv(1.6.1),
one gets 
$$
\E[\t_{I}^{m}]=
\E[\t_{I}^{m}|\t_{I}^m<\t_{m}^{m}]+
{\E[\t_{m}^{m},\t_{m}^{m}=\t_{I}^{m}]
\over
\P[\t_{I}^{m}<\t_{m}^{m}]}
\Eq(2.6.3)
$$
Bounding the first summand on the right by Theorem \thv(LL.7)
gives \eqv(2.6.1). Using Theorem \thv(LL.5) for the right hand side of 
\eqv(2.6.1) gives \eqv(2.6.2). \endproof

\newpage

\chap{4. Laplace transforms and spectra}4

In this section we present a characterization of the spectrum 
of the Dirichlet operator $(1-P_{N})^{I}$,    
$I\sb\MM_{N}$,  in terms of 
Laplace transforms of transition times (defined in 
\eqv(1.3.3) and \eqv(1.3.4)). This  connection forms the basis 
of the investigation of the low-lying 
spectrum that is presented in Section 5. To exploit this 
characterization we study the region of 
analyticity and boundedness of Laplace transforms. As a 
first consequence we then 
show that the  principal eigenvalue 
for Dirichlet operators 
are with high precision equal to the inverse of expected 
transition times. A combination  of these results 
then leads to the characterization of the 
low-lying spectrum given in the next section. 

For any $J\sb\MM_{N}$ we denote the 
principal eigenvalue of the Dirichlet-operator $P_N^J$ by
$$
\l_{J}\equiv\min\s((1-P_{N})^{J})
\Eq(3.1.1)
$$
For $I,J\sb\MM_N$ we define the matrix 
$$
\GG_{I,J}(u)\equiv \left(
\d_{m^\prime,m}-G_{m,I\cup J}^{m^\prime}(u)
\right)_{m^\prime,m\in J\ba I}
\Eq(3.1.1a)
$$
where $\d_{x,y}$ is Kronecker's symbol. 
We then have

\lemma{\TH(C.1)} {\it
Fix subsets $I,J\sb\MM_{N}$ such that $J\ba I\neq\em$ 
and a number 
$0\leq\l\equiv 1-e^{-u}<\l_{I\cup J}$. Then 
$$
\l\in\s((1-P_{N})^{I})
\qquad\Longleftrightarrow\qquad
\det\GG_{I,J}(u)=0
\Eq(3.1.2)
$$
Moreover, the map 
$\ker\GG_{I,J}(u)\ni\vec{\phi}\mapsto
\phi\in\1_{I^c}\R^{\G_{N}}$ defined by 
$$
\phi(x)\equiv 
\sum_{m\in J\ba I} \vec{\phi}_m K_{m,I\cup J}^x(u),
\qquad x\in\G_{N}
\Eq(3.1.3)
$$
is an isomorphism onto the eigenspace corresponding to 
the eigenvalue $\l$.
}

\proof
Assume that $\phi$ is an eigenfunction with 
corresponding eigenvalue $\l<\l_{I\cup J}$. 
We have to prove that $\GG_{I,J}(u)$ is 
singular. In view of  
\eqv(1.3.4a) the condition $\l_{I\cup J}>\l$ 
implies that $\tilde{\phi}$ defined below is finite.
$$
\tilde{\phi}\equiv
\sum_{m\in J} \phi(m) K_{m,I\cup J}^x(u),
\qquad x\in\G_{N}
\Eq(3.1.3a)
$$ 
Furthermore, \eqv(1.5.1) and \eqv(1.3.4) 
imply for $x\in\G_{N}$
$$
e^{u}(1-P_{N}-(1-e^{-u}))
\tilde{\phi}(x)=
(1-e^{u}P_{N})\tilde{\phi}(x)=
\sum_{m^\prime\in I\cup J}\d_{m^\prime,x}
\sum_{m\in J}\phi(m)
\left(
\d_{m^\prime,m}-G_{m,I\cup J}^{m^\prime}(u)
\right)
\Eq(3.1.4)
$$
Let $\D\equiv \phi-\tilde{\phi}$. We want to show 
$\D=0$. Obviously, we have $\D$ vanishes on $I\cup J$ and 
$\tilde{\phi}$ on $I$. Combining \eqv(3.1.4) 
with the eigenvalue equation for $\phi$ and 
the choice of $u$, we obtain
$$
\eqalign{
(1-P_{N})^{I\cup J}\D&=
\1_{(I\cup J)^c}(1-P_{N})^{I}\D=\1_{(I\cup J)^c}
\left(
(1-P_{N})^{I}\phi-(1-P_{N})\tilde{\phi}
\right)
\cr &=
\1_{(I\cup J)^c}
(\l\phi-(1-e^{-u})\tilde{\phi})=\l\D
}
\Eq(3.1.5)
$$
Since $\l\notin\s((1-P_{N})^{I\cup J})$, we 
conclude $\D=0$. 
Replacing $\tilde{\phi}$ by $\phi$ in \eqv(3.1.4) 
and, using $\l\equiv 1-e^{-u}$ again, gives
$$
0=
\sum_{m^\prime\in I\cup J}\d_{m^\prime,x}
\sum_{m\in J}\phi(m)
\left(
\d_{m^\prime,m}-G_{m,I\cup J}^{m^\prime}(u)
\right),
\qquad x\in I^c
\Eq(3.1.6)
$$
Choosing $x\in J\ba I$ yields that
$(\phi(m))_{m\in J\ba I}\in\ker\GG_{I,J}(u)$ 
and the right-hand side of the equivalence 
in \eqv(3.1.2) follows. In particular, we have proven 
that the restriction map 
$\phi\mapsto(\phi(m))_{m\in J\ba I}$ defined on the 
eigenspace corresponding to $\l$ is the inverse of 
the map defined in \eqv(3.1.3). 

For the converse implication we note that 
for  $\l<\l_{I\cup J}$ the entries of the matrix 
$\GG_{I,J}(u)$ are finite. We replace 
$(\phi(m))_{m\in J\ba I}$ in \eqv(3.1.3a) by the 
solution $\vec{\phi}$ of the linear system 
$\GG_{I,J}(u)\vec{\phi}=0$ and deduce from 
\eqv(3.1.4) and \eqv(3.1.6) that $\l$ is an 
eigenvalue with eigenfunction $\tilde{\phi}$.
\endproof

As a first step we now derive a  lower bound 
on these eigenvalues, using a Donsker-Varadhan [DV] like 
argument that we will later prove to be sharp. 

\lemma{\TH(C.2)} {\it
For every nonempty subset $J\sb\MM_{N}$ we have
$$
\l_{J}\max_{x\notin J}\E[\t_J^x]\geq 1
\Eq(3.2.1)
$$
}

\proof
For $\phi\in\R^{\G_{N}}$ we have for all 
$x,y\in\G_N$ and $C>0$ 
$$
\phi(y)\phi(x)\leq
{1\over 2}(\phi(x)^2C+\phi(y)^2/C)
\Eq(3.2.2)
$$
Thus choosing $C\equiv \psi(y)/\psi(x)$, 
where $\psi\in\R^{\G_{N}}$ is such that 
$\psi(x)>0$ for all $x\in\op{supp}\phi$,
we compute, using reversibility,
$$
\eqalign{
\la P_{N}\phi,\phi\ra_{\Q_{N}}
\leq&{1\over 2}\sum_{x,y\in\G_{N}}
\Q_{N}(x)P_{N}(x,y)(\phi(x)^2(\psi(y)/\psi(x))+
\phi(y)^2(\psi(x)/\psi(y)))
\cr=& \sum_{x,y\in\G_{N}}
\Q_{N}(x)\phi(x)^2
{P_{N}(x,y)\psi(y)\over \psi(x)}
= \left\la \phi\biggl(
{P_{N}\psi\over \psi}\biggr),
\phi\right\ra_{\Q_{N}}
}
\Eq(3.2.3)
$$
Let $\phi$ be an eigenfunction for the principal 
eigenvalue and set $\psi(x)\equiv \E[\s_{J}^x]$, 
$x\in\G_{N}$. Invoking \eqv(1.5.2) for $u\equiv 0$ 
and $I\equiv J$ we get
$$
\l_{J}||\phi||_{\Q_{N}}^2\geq
\la\phi/\psi,\phi\ra_{\Q_{N}}
\Eq(3.2.4)
$$
which in turn gives the assertion.
\endproof

We now 
study the behavior of Laplace transforms slightly 
away from their first pole on the real axis.

\lemma{\TH(C.4)} {\it
Fix nonempty subsets $I,J\sb\MM_{N}$. Let 
$G_{I,J}^{x}$ be the Laplace transform defined in 
\eqv(1.3.3). It follows that for some $c>0$ and for 
$k=0,1$ uniformly in
$0\leq \Re(u),|\Im(u)|\leq c/(c_N T_{I\cup J})$ 
and $x\in\G_N$
$$
\partial_u^k{G}_{I,J}^x(u)=
\left(
1+\OO(|u|c_NT_{I\cup J})
\right)
\partial_u^k{G}_{I,J}^x(0)
\Eq(3.4.1)
$$
}

\proof
By \eqv(1.3.4a), we know that 
$G_{I,J}^x(u)$, 
$x\in\G_{N}$, are finite for all $u$ such that 
$1-e^{-\Re(u)}<\l_{I\cup J}$. Put 
$$
K_{u,v}\equiv
K_{I,J}^{(\cdot)}(u)-K_{I,J}^{(\cdot)}(v)
\Eq(3.4.2)
$$
\eqv(1.5.1) and \eqv(1.5.2) imply that for $k=0,1$,
$$
(1-P_{N})^{I\cup J}
(\partial_{u}\partial_{v})^k{K}_{u,0}=
(1-e^{-u})\partial_u^kK_{I,J}^{(\cdot)}(u)+
\d_{k,1}K_{u,0}
\Eq(3.4.3)
$$
We first 
consider the case where $k=0$. Using \eqv(2.2.5), we get from 
\eqv(3.4.3)  for all $x\notin I\cup J$
$$
{G_{I,J}^{x}(u)\over G_{I,J}^{x}(0)}
=
1+(1-e^{-u})
\sum_{y\notin I\cup J}
{\P[\s_{y}^{x}<\t_{I}^{x}|\t_{I}^{x}\leq\t_{J}^{x}]
\over
\P[\t_{I\cup J}<\t_{y}^{y}]}
\D_{I\cup J}(x,y)
{G_{I,J}^{y}(u)\over G_{I,J}^{y}(0)}
\Eq(3.4.5)
$$
where $\D_{I\cup J}$ is defined in  \eqv(2.2.6).
Setting 
$$
M_{N,k}(u)\equiv\max_{x\notin I\cup J}
{|\partial_u^kG_{I,J}^{x}(u)|\over G_{I,J}^{x}(0)}
\Eq(3.4.6)
$$
and, using that $\frac{\partial_u^kG_{I,J}^{x}(u)}{ G_{I,J}^{x}(0)}
=\E[\t^x_I|\t^x_I<\t^x_J]$, we obtain from \eqv(3.4.5) that for
$1-e^{-\Re(u)}<\l_{I\cup J}$
$$
1-|1-e^{-u}|M_{N,0}(u)M_{N,1}(0)
\leq M_{N,0}(u)\leq
1+|1-e^{-u}|M_{N,0}(u)M_{N,1}(0)
\Eq(3.4.7)
$$
But by Theorem \thv(LL.7) we have a uniform bound on $M_{N,1}(0)$, 
and this implies \eqv(3.4.1) for $x\not\in I\cup J$.

For $k=1$ \eqv(3.4.3) gives
$$
{\del_uG_{I,J}^{x}(u)\over G_{I,J}^{x}(0)}
=
{\del_uG_{I,J}^{x}(0)\over G_{I,J}^{x}(0)}
+
\sum_{y\notin I\cup J}
{\P[\s_{y}^{x}<\t_{I}^{x}|\t_{I}^{x}\leq\t_{J}^{x}]
\over
\P[\t_{I\cup J}<\t_{y}^{y}]}\D_{I\cup J}(x,y)
\biggl(
(1-e^{-u}){\del_u {G}_{I,J}^{y}(u)\over G_{I,J}^{y}(0)}
+{G_{I,J}^{y}(u)\over G_{I,J}^{y}(0)}-1
\biggr)
\Eq(3.4.5a)
$$
and the same arguments together  with \eqv(3.4.1) 
for $k=0$ show, for some $c>0$ and all 
$0\leq\Re(u),|\Im(u)|<c c^{-1}T^{-1}_{J\cup I}$, that 
$$
M_{N,1}(u)\leq
M_{N,1}(0)\left(1+\OO(|u|c_NT_{I\cup J})\right)+
|1-e^{-u}|M_{N,1}(u)M_{N,1}(0)
\Eq(3.4.7a)
$$
In particular, we conclude that on the same set, 
$$
M_{N,1}(u)=\OO(M_{N,1}(0))=
\OO(c_NT_{I\cup J})
\Eq(3.4.8)
$$
Inserting this estimate into \eqv(3.4.5a) 
\eqv(2.3.1) and \eqv(3.4.1) for $k=0$ again gives for 
all $0\leq \Re(u),|\Im(u)|<c c_N T_{I\cup J}$
$$
{\del_u {G}_{I,J}^{x}(u)\over G_{I,J}^{x}(0)}=
\left(
1+\OO(|u|Ne^{Nd_{J\cup K}})
\right)
{\del_u {G}_{I,J}^{x}(0)\over G_{I,J}^{x}(0)},
\qquad x\notin I\cup J
\Eq(3.4.9)
$$
which yields \eqv(3.4.1) for $k=1$ and 
$x\notin I\cup J$. 

The remaining part, namely $x\in I\cup J$, 
follows by first using  \eqv(1.5.1), respectively \eqv(1.5.2),
to express the quantities $\del^kG_{I,J}^x$ in terms 
of $\del^kG_{I,J}^y$ with $y\not\in I\cup J$ and then applying the 
result obtained before.
\endproof

We now have all tools to establish a sharp relation between mean exit times and
the principal eigenvalue  $\l_I$ of $P_N^I$. 
Set $u_I\equiv -\ln(1-\l_I)$.  
We want to show that  
$$
G_{m,I}^{m}(u_{I})=1
\Eq(3.5.1)
$$
Indeed, this follows from Lemma \thv(C.1) with $J=I\cup \{m\}$, $m\in\MM_N$,
if we can show that $\l_I<\l_{I\cup m}$. Now it is obvious by monotonicity that
$\l_I\leq\l_{I\cup m}$. But if equality held, then by \eqv(1.3.4a),
$\lim_{u\uparrow u_I} G_{m,I}^m(u)=+\infty$; by continuity, it follows that
there exists $u<u_I$ such that $G_{m,I}^m(u)=1$, implying by Lemma \thv(C.1)
that $1-e^{-u}<\l_I$ is an eigenvalue of $P_N^I$, contradicting the fact that
$\l_I$ is the smallest eigenvalue of  $P_N^I$. We must conclude that
$\l_I<\l_{I\cup m}$ and that \thv(3.5.1) holds.

\theo{\TH(C.5)}{\it
Fix a proper nonempty subset $I\sb\MM_{N}$. 
Let $m\in\MM_{N}\ba I$ be the unique 
local minimum satisfying 
$
T_{I}=T_{m,I}$.
Then 
$$
\l_{I}=\left(
1+\OO(T_{I\cup m}/T_{I})
\right)
\E[\t_{I}^{m}]^{-1}
\Eq(3.5.3)
$$
In particular, 
$$
\l_{I}= R_m T_{I}^{-1}\left(1+\OO(\e_N|\G_N|+
|\G_N|/(\e_N\a_Nb_N))\right)
\Eq(3.5.4)
$$
}

\proof 
Using that for $x\geq 0$, $e^x>1+x$, for real and positive $u$,
$$
G^m_{m,I}(u)=\E\left[ e^{u\t^m_m}\1_{\t^m_m<\t^m_I}\right]
\geq \P[\t^m_m<\t^m_I] +u \E\left[ \t^m_m\1_{\t^m_m<\t^m_I}\right]
\Eq(3.5.0)
$$
Using this in \eqv(3.5.1), we 
immediately obtain the upper bound 
$$
u_{I}\leq 
\frac{\P[\t_{I}^{m}<\t_{m}^{m}]}
{\E\left[\t_{m}^{m}\1_{\t_{m}^{m}<\t_{I}^{m}}\right]}
\Eq(3.5.5)
$$
Using now Lemma \thv(LL.8) to bound the right hand side, gives 
the upper bound of \eqv(3.5.3). The lower bound is of course already 
contained in Lemma \thv(C.2).
\endproof

The a priori control of the Laplace transforms 
given in Lemma \thv(C.4) can be used to control 
denominators in the renewal relation \eqv(1.6.1) 
which will be important for the construction 
of the solution of the equation appearing in 
\eqv(3.1.2). We are interested in the 
behavior of $G_{m,I}^{m}$ near $u_{I}$. 

\lemma{\TH(C.7)} {\it
Under the hypothesis of  Theorem \thv(C.5) there exists 
$c>0$ such that for all 
$0\leq \Re(u)<c/( c_NT_{I\cup m})$
$$\eqalign{
G_{m,I}^m(u)-1&=\E\left[\t^m_m\1_{\t^m_m<\t^m_I}\right]
\left(
u-u_{I}+(u-u_{I})^{2}\OO(c_NT_{I\cup m})
\right)\cr
&=(1+\OO(\e_N)) R_m^{-1} 
\left(
u-u_{I}+(u-u_{I})^{2}\OO(c_NT_{I\cup m})
\right)
}
\Eq(3.7.1)
$$
}

\proof
Performing a Taylor expansion at $u=u_{I}$ to 
second order of the Laplace transform on the 
left-hand side of \eqv(3.7.1) and recalling 
\eqv(3.5.1) we get
$$
G_{m,I}^{m}(u)-1=\del_u {G}_{m,I}^m(u_{I})
\left((u-u_{I})
-(u-u_{I})^2\RR_{I}(u)\del_u {G}_{m,I}^m(u_{I})^{-1}
\right)
\Eq(3.7.2)
$$
where
$$
\RR_{I}(u)\equiv\int_{0}^{1}s
\ddot{G}_{m,I}^{m}((1-s)u_{I}+su)ds
\Eq(3.7.3)
$$
\eqv(3.7.2) then follows from Cauchy's 
inequality combined with \eqv(3.4.1) 
and \eqv(3.5.4) which shows, 
for $c>0$ small enough, $C<\infty$ large 
enough, and all $u$ considered in the Theorem, that
$$
\eqalign{
|\ddot{G}_{m,I}^m(u)|\leq
\ddot{G}_{m,I}^m\left(c /(c_NT_{I\cup m})\right)
\leq & C c_NT_{I\cup m}
\del_u {G}_{m,I}^{m}\left(c /(c_NT_{I\cup m})\right)
\cr \leq & 
C^2 c_NT_{I\cup m}\del_u {G}_{m,I}^m(0)
}
\Eq(3.7.4)
$$
where we used Lemma \thv(C.4). Using Lemma \thv(LL.8), the assertion of the
lemma follows. 
\endproof

\newpage

\chap{5. Low lying eigenvalues}5

In the present section we prove the main new result of this paper. Namely,
we establish a precise relation between the low-lying part of the spectrum of 
the
operator $1-P_N$ and the metastable exit times associated to the set $\MM_N$.
Together with the results of Section 2, this allows us to give sharp
estimates on the entire low-lying spectrum in terms of the transition 
probabilities  between points in $\MM_N$ and the invariant measure.

As a matter of fact we will prove a somewhat more general result. Namely, 
instead of computing just the low-lying spectrum of $1-P_N$, we will 
do so for any of the Dirichlet operators
$(1-P_N)^I$, with $I\subset \MM_N$ (including the case $I=\emptyset$).
In the sequel we will fix 
$I\sb\MM_{N}$ with  $I\neq \MM_{N}$. 

The strategy of our proof will be to show that to each of the points
$m_i\in\MM_N\ba I$ corresponds exactly one eigenvalue $\l^I_i$ 
of $(1-P_N)^I$ and that
this eigenvalue in turn is close to the principle eigenvalue of some
Dirichlet operator $(1-P_N)^{\S_i}$, with $I\subset\S_i\subset\MM_N$.
We will now show how to construct these sets $\S_i$ in such a way as to obtain 
an ordered sequence of eigenvalues. 

We set the first 
exclusion set $\S_{0}$ and the first effective depth $T_{1}$ 
to be 
$$
\S_{0}\equiv I \text{and} T_{1}\equiv T_{\S_{0}}
\Eq(4.1.1f)
$$
where $T_{K}$, $K\sb\MM_{N}$, 
is defined in \eqv(0.61). If $I\neq\emptyset$, let  $m_{1}$ 
be the unique  point in $\MM_N\ba I$ such that 
$$
T_{m_{1},I}=T_{1}
\Eq(4.1.1a)
$$
If $I= \em$, let  $m_1$ be the unique 
element of $\MM_N$ such that $\Q_N(m_1)=\max_{m\in\MM_N}\Q_N(m)$.

 For $j=2,\ldots,j_{0}$, 
$j_{0}\equiv |\MM_{N}\ba I|$, we define the corresponding quantities 
inductively by 
$$
\S_{j-1}\equiv \S_{j-2}\cup m_{j-1} \text{and} 
T_{j}\equiv T_{\S_{j-1}}
\Eq(4.1.1d)
$$
and 
$m_{j}\in\MM_{N}\ba \S_{j-1}$ is 
determined by the equation
$$
T_{N}(m_{j},\S_{j-1})=T_{j}
\Eq(4.1.1e)
$$
In order to avoid distinction as to whether or not 
$j=j_{0}$, it will be convenient to set 
$T_{j_{0}+1}\equiv b_N^{-1}$.
Note that this construction and hence all the sets $\S_j$ depend on 
$N$.  
An important fact is that the sequence $T_j$ is decreasing.
To see this, note that by construction and the assumption of genericity
$$
T_l=T_{m_l,\S_{l-1}}\geq \e_N^{-1}T_{m_{l+1},\S_{l-1}}\geq
\e_N^{-1} T_{m_{l+1},\S_l}=\e_N^{-1}T_{l+1}
\Eq(UU.0)
$$
The basic  heuristic picture behind this construction 
can be summarized as follows. 
To each $j=1,\ldots,j_0$ associate a  rank one operator 
obtained by projecting the Dirichlet operator 
$(1-P_N)^{\S_{j-1}}$ onto the eigenspace corresponding 
to its principal eigenvalue $\l_{\S_{j-1}}\sim T_j^{-1}$. 
Note that our construction of $\S_j$ as an increasing sequence 
automatically guarantees that these eigenvalues will be in increasing order.
The direct 
sum of these  rank one operators acts approximately like 
$(1-P_{N})^{I}$ on the eigenspace corresponding to 
the exponentially small part of its spectrum. Hence 
the difference between both operators can be treated 
as a small perturbation. 

\remark We can now explain what the minimal non-degeneracy conditions
are that are necessary for Theorem \thv(A.1) to hold. Namely, what must be 
ensured is that the preceding construction of the sequence of sets 
is {\it unique}, and that the $T_{\S_j}$ are by a diverging factor
$\e^{-1}_N$ larger than all other $T_{x,\S_j}$. 

We are now ready to formulate the main theorem of this section.
Let $\l_{j}$, $j=1,\ldots,|\G_{N}\ba I|$, be the 
$j$-th eigenvalue of $(1-P_{N})^{I}$ written in 
increasing order and counted with multiplicity 
and pick a corresponding eigenfunction $\phi_{j}$ 
such that $(\phi_j)_j$ is an orthonormal basis of 
$\1_{I^c}\ell^{2}(\G_N,\Q_N)$. We 
then have 

\theo{\TH(S.1)} {\it
Set $j_0\equiv|\MM_N\ba I|$. There is $c>0$ such 
that the Dirichlet operator $(1-P_{N})^{I}$ 
has precisely $j_{0}$ simple eigenvalues in the 
interval $[0,cb_N)|\G_N|$, i.e.
$$
\s((1-P_{N})^{I})\cap [0,c b_N|\G_N|^{-1})
=\{\l_{1},\ldots,\l_{j_{0}}\}
\Eq(4.1.2)
$$
Define $\TT_1\equiv\infty$ and for $j=2,\ldots, j_{0}$
$$
\TT_{j}\equiv\min_{1\leq k<j} 
T_{m_{k},m_{j}}/T_{j}\geq \e_N^{-1}
\Eq(4.1.3a)
$$
Then
$$
\l_j=\left(
1+\OO(\TT_{j}^{-1}+T_{j+1}/T_{j}))
\right)
\l_{\S_{j-1}}
\Eq(4.1.4)
$$
where $\l_K$, $K\sb\MM_N$, is defined in \eqv(3.1.1). 

Moreover, the eigenfunction $\phi_{j}$ satisfies 
for $k=1,\ldots,j-1$
$$
\phi_{j}(m_{k})=\phi_{j}(m_{j})
\OO\left(R_{m_j}T_{m_{k},m_{j}}/T_{j}\right)
\Eq(4.1.6)
$$
}

\remark Combining Theorem \thv(S.1) with Theorem \thv(C.5) 
and Theorem \thv(LL.5), we get immediately 

\corollary{\TH(S.1a)}{\it With the notation of Theorem \thv(S.1),
for $j=1,\ldots,j_0$ that
$$
\eqalign{
\l_{j}&=\left(1+
\OO(\TT_{j}+T_{j+1}/T_{j})
\right)
\E\left[\t_{\S_{j-1}}^{m_{j}}\right]^{-1}
\cr 
&= \frac 1{T_{j}} R_{m_j} \left(1+\OO\left(
|\G_N|(\e_N+1/(a_N b_N\e_N))\right)\right)
}
\Eq(4.1.3)
$$
}
Note that Corollary \thv(S.1a) is a precise version of (ii) of Theorem 
\thv(A.1). The estimate \eqv(4.1.6), together with 
the representation \eqv(3.1.3) and the estimates of the Laplace transforms
in Lemma \thv(C.4), gives a precise control of the eigenfunctions
and implies in particular (iv) of Theorem \thv(C.4).

The strategy of the proof will be to seek, for 
each $J\equiv \S_{j}$, for a solution of the equation 
appearing in \eqv(3.1.2) with $\l$ near the principle eigenvalue
 of the associated Dirichlet operator $(1-P_{N})^{\S_{j-1}}$. 
We then show that these eigenvalues are simple and that 
no other small eigenvalues occur. 

For the investigation of the structure of the equations 
written in \eqv(3.1.2) we have to take a closer look at 
the properties of the effective depths defined in 
\eqv(4.1.1d). We introduce for all $m\in\MM_N\ba I$ the 
associated ``metastable depth'' with exclusion at $I$ by 
$$
T_N(m)\equiv T_{m,\MM_N(m)}, \text{where}
\MM_N(m)\equiv 
I\cup\{m^\prime\in\MM_N\,|\,\Q_N(m^\prime)>\Q_N(m)\}
\Eq(meta1)
$$ 
Let us define for $j=2,\ldots,j_0$
$$
\EE_{j}\equiv \min_{1\leq l<j}
T_{m_{l},\S_{j}\ba m_l}
\Eq(4.2.1)
$$

The following result relates our inductive definition to these
geometrically more transparent objects and establishes some crucial 
properties:

\lemma{\TH(S.2)} {\it
Every effective depth is a metastable depth, more precisely 
for all $j=1,\ldots,j_0$ it follows
$$
T_j=T_N(m_j)(1+\OO(\e_N|\MM_N|))
\Eq(meta2)
$$
For $j=2,\ldots,j_0$ we have 
$$
\TT_{j}\geq \EE_{j}/T_{j}\geq\e_N^{-1}.
\Eq(4.2.2)
$$ 
Moreover, for $j,l=1,\ldots,j_0$, $l<j$, we have 
$$
T_{m_{l},\S_{j}\ba m_{l}}=
T_{\S_{j}\ba m_{l}}(1+\OO(\e_N|\MM_N|)) 
\Eq(4.2.2a)
$$
}

\proof 
Fix $l<j$. It will be convenient to decompose 
$\S_j=\S_{l-1}\cup m_l\cup \S_j^+$, where
$\S_j^+\equiv \S_j\ba \S_l$.  We will use heavily the
(almost) ultra-metric
$e(\cdot,\cdot)$ introduced in Section 2; for the purposes of the proof 
we can ignore the irrelevant errors in the ultra-metric inequalities
(i.e. all equalities and inequalities relating the functions
$e$ in the course of the proof are understood up to error of
at most $\ln 3$).
Note that $\ln T_{x,J} =e(x,J)-f(x)$, where $f(x)\equiv -\ln \Q_N(x)$. 
In particular, $d_l\equiv\ln T_l=e(m_l,\S_{l-1})-f(m_l)$.
As a first step we prove the following general fact that will be used several 
times:

\lemma{\TH(UU.1)}{\it Let $m$ be such that $e(m,m_l)< e(m_l,\S_{l-1})$. 
Then 
$f(m)\geq f(m_l)+|\ln \e_N|$.}

\proof Note that by ultra-metricity,
$$
e(m,\S_{l-1})=\max\left(e(m,m_l),e(m_l,\S_{l-1})\right)
=e(m_l,\S_{l-1})
\Eq(UU.001)
$$
But since for any $m$,
$$
e(m,\S_{l-1})-f(m)\leq d_l-|\ln \e_N|=e(m_l,\S_{l-1})-f(m_l)-|\ln \e_N|
\Eq(UM.2)
$$
which implies by \eqv(UU.001)  $f(m_l)\leq f(m)-|\ln \e_N|$.\endproof

Let us now start by proving \eqv(4.2.2). The first inequality is trivial.
We distinguish the cases where $e(m_l,\S_j^+)$ is larger or smaller than 
$e(m_l,\S_{l-1})$.
\item{(i)} Let 
$e(m_l,\S_j^+)\geq e(m_l,\S_{l-1})$. \hfill\break
Since $e(m_l,\S_j\ba m_l)=\min\left(e(m_l,\S_{l-1}),e(m_l,\S_j^+)\right)$,
this implies that $e(m_l,\S_j\ba m_l)=e(m_l,\S_{l-1})$. 

Then, using \eqv(UU.0) and genericity from Definition 1.2, 
$$
\eqalign{
e(m_l,\S_j\ba m_l)-f(m_l)&=e(m_l,\S_{l-1})-f(m_l)=d_l
\geq e(m_{j-1},\S_{l-1})-f(m_{j-1})\cr
&\geq e(m_{j-1}\S_{j-2})
-f(m_{j-2})=d_{j-1}\geq d_j+|\ln \e_N|
}
\Eq(UM.1)
$$
Obviously, this gives \eqv(4.2.2) in this case. 
\item{(ii)} 
 Let 
$e(m_l,\S_j^+)< e(m_l,\S_{l-1})$. \hfill\break 
In this case there must exist $m_k\in\S_j^+$ such that
$e(m_l,\S_j\ba m_l)=e(m_l,m_k)$, and hence
$e(m_k,m_l)<e(m_l,\S_{l-1})$. 
Thus we can use Lemma \thv(UU.1) for 
 $m=m_k$.  Together with the trivial
inequality $e(m_k,m_l)\geq e(m_k,\S_{k-1})$, it follows that
$$
\eqalign{
&e(m_l,\S_{j}\ba m_l)-f(m_l)=e(m_k,m_l)-f(m_l)\cr
&\geq 
e(m_k,\S_{k-1})-f(m_k)+f(m_l)-f(m_k)\geq d_k+|\ln \e_N|\geq
d_j+|\ln \e_N|
}
\Eq(UM.3)
$$
This implies \eqv(4.2.2) in that case and concludes the proof
of this inequality.

We now turn to the proof of \eqv(4.2.2a). We want to proof
that the maximum over $T_{m,\S_j\ba m_l}$ is realized for $m=m_l$.
Note first that it is clear that the maximum cannot be realized for $m\in
\S_j\ba m_l$ (since in that case $T_{m,\S_j\ba m_l}=1$). Thus fix 
$m\not\in \S_j$. We distinguish the cases $e(m,m_l)$ less or larger than 
$e(m,\S_{j}\ba m_l)$. 
\item{(i)} Assume $e(m,m_l)<e(m,\S_j\ba m_l)$.\hfill\break
The ultra-metric property of $e$ then implies that
$e(m_l,\S_j\ba m_l)=e(m,\S_j\ba m_l)$, and hence, using the 
 argument from above,
$f(m)>f(m_l)+|\ln \e_N|$. Thus
$$
e(m_l,\S_j\ba m_l)-f(m_l)=e(m,\S_j\ba m_l)-f(m)+f(m)-f(m_l)
\geq e(m,\S_j\ba m_l)-f(m) +|\ln \e_N|
\Eq(UM.5)
$$
which excludes that in this case $m$ may realize the maximum.
We turn to the next case. 
\item{(ii)} Assume  $e(m,m_l)\geq e(m,\S_j\ba m_l)$.\hfill\break
We have to distinguish the two sub-cases like in the proof of \eqv(4.2.2).
\itemitem{(ii.1)} $e(m_l,\S_j^+)\geq  e(m_l,\S_{l-1})$. \hfill\break 
Here we note simply that by \eqv(UM.1) 
$$
e(m_l,\S_j\ba m_l)=e(m_l,\S_{l-1})-f(m_l)=d_l>e(m,\S_{l-1})-f(m)
\geq e(m,\S_j\ba m_l)-f(m)
\Eq(UM.6)
$$
which implies that $m$ cannot be the maximizer.
\itemitem{(ii.2)} $e(m_l,\S_j^+) <  e(m_l,\S_{l-1})$. \hfill\break 
This time we use  \eqv(UM.3)  for some $m_k\in \S_j^+$ and so
$$
e(m_l,\S_j\ba m_l)-f(m_l)>d_k>e(m,\S_{k-1})-f(m)
 \geq e(m,\S_j\ba m_l)-f(m)
\Eq(UM.7)
$$
where in the last inequality we used that by assumption 
$e(m,m_l) >e(m,\S_j\ba m_l)$. Again \eqv(UM.7) rules out $m$ as maximizer, and
since all cases are exhausted, we must conclude that \eqv(4.2.2a) holds.

It remains to show that \thv(meta2) holds. Now the crucial observation is that
by Lemma \thv(UU.1), 
$$
\MM_N(m_j)\cap \left\{m\in\MM_N: e(m_j,m)<e(m_j,\S_{j-1})\right\}=\emptyset
\Eq(meta3)
$$
Thus, for all $m\in\MM_N(m_j)$, $T_{m_j,m}\geq T_{m_j,\S_{j-1}}$, which implies 
of course that 
$$
T_{m_j,\MM(m_j)}\geq 
 T_{m_j,\S_{j-1}}
\Eq(UM.8)
$$
To show that the converse inequality also holds,
it is obviously enough to show that the set
$$
\{m|T_{m_j,m}\leq T_{m_j,\S_{j-1}}\}\cap \MM_N(m_j)\neq\emptyset
\Eq(UM.9)
$$
Assume the contrary, i.e. that for all $m\in\MM(m_j)$ $T_{m_j,m}>
T_{m_j,\S_{j-1}}$. 
Now let $m\not\in I$ be such a point. Then also 
$e(m_j,m)>e(m_j,\S_{j-1})$, and so by ultra-metricity
$e(m,\S_{j-1})=\max\left(e(m_j,m),e(m_j,\S_{j-1})\right)
>e(m_j,\S_{j-1})$. 
But, since $f(m)\leq f(m_j)$,
it follows that
$$
T_{m,\S_{j-1}}>T_{m_j,\S_{j-1}}
\Eq(UM.10)
$$
in contradiction with the defining property of $m_j$. Thus \eqv(UM.9)
must hold, and so $T_{m_j,\MM_N(m_j)}\leq   T_{m_j,\S_{j-1}}$.
This concludes the proof of the Lemma.\endproof

We now turn to the constructive part of the investigation 
of the low lying spectrum. 
Having in mind the heuristic picture described before 
Theorem \thv(S.1) we are searching for solutions $u$ of 
\eqv(3.1.2) for $J\equiv \S_{j}$ near $u_{\S_{j-1}}\equiv 
-\log(1-\l_{\S_{j-1}})$. The procedure of finding $u$ 
is as follows. The case $j=1$ 
was studied in Theorem 3.5. For $j=2,\ldots,j_{0}$ we 
consider the matrices $\GG_{j}=\GG_{I,\S_{j}}$ 
defined in \eqv(3.1.1a), i.e.
$$
\GG_{j}\equiv\left(
\matrix\format\c\quad&\c\\
\KK_{j} & -\vec{g}_{j}\\
-(\vec{g}_{j})^t & 1-G_{m_{j},\S_{j}}^{m_{j}}
\endmatrix
\right)
\equiv
\left(
\matrix\format \c&\c&\c&\c\\
1-G_{m_1,\S_j}^{m_1} &  -G_{m_2,\S_j}^{m_1} &          \hdots           &   -G_{m_j,\S_j}^{m_1}  \\
-G_{m_1,\S_j}^{m_2}  &      \ddots          &                           &         \vdots         \\
       \vdots        &                      &                           &                        \\
                     &                      &                           &-G_{m_j,\S_j}^{m_{j-1}} \\
-G_{m_1,\S_j}^{m_j}  &      \hdots          & -G_{m_{j-1},\S_j}^{m_{j}} &  1-G_{m_j,\S_j}^{m_j}  \\
\endmatrix
\right)
\Eq(4.3.1)
$$
and define 
$$
\NN_{j}\equiv \DD_{j}-\KK_{j},\text{where} 
\DD_{j}\equiv\op{diag}(1-G_{m_l,\S_j}^{m_l})_{1\leq l<j}
\Eq(4.3.2)
$$

Equipped with the structure of the effective depths 
written in Lemma \thv(S.2) and the control of Laplace 
transforms of transition times obtained in the previous 
chapter one simply can write a Neumann series for 
$\1-\DD_{j}(u)^{-1}\NN_{j}(u)$ for $u$ near $u_{\S_{j-1}}$ 
proving 
the invertibility of $\KK_{j}(u)$. We then compute  
$$
\det\GG_{j}=
\det
\pmatrix\format\c&\c\\
\KK_j             & 0       \\
-(\vec{g}_j)^t\quad & G_{j}
\endpmatrix
=G_{j}\det\KK_{j}
\Eq(4.3.2d)
$$
where 
$$
G_{j}\equiv 1-G_{m_{j},\S_{j}}^{m_{j}}-
(\vec{g}_{j})^t\KK_{j}^{-1}\vec{g}_{j}
\Eq(4.3.2e)
$$
This follows  by simply adding the column vector
$$
\pmatrix
\KK_j\\-(\vec{g}_j)^t
\endpmatrix
\KK_j^{-1}\vec{g}_j
$$
(which clearly is a linear combination of the first $j-1$
columns of $\GG_j$) to the last column in $\GG_j$, and the fact that this 
operation leaves 
the determinant unchanged.
From this representation we construct 
solutions $\tilde{u}_{j}$ near $u_{\S_{j-1}}$ of \eqv(3.1.2). 
We begin 
with

\lemma{\TH(S.3)} {\it
For all $j=2,\ldots,j_{0}$ there 
are constants $c>0$, $C<\infty$ such that for all 
$C^\prime<\infty$ and all 
$$
CR_{m_j}\EE^{-1}_{j}<\Re(u)<c c^{-1}_NT^{-1}_{j+1},
\qquad
|\Im(u)|<c/(c_N T_{j+1})
\Eq(4.3.3)
$$
the inverse of $\KK_{j}(u)$ exists. 
The $l$-th component of $\KK_{j}(u)^{-1}\vec{g}_{j}(u)$ 
restricted to 
the real axis is strictly monotone increasing and,  
uniformly in $u$,
$$
(\KK_{j}(u)^{-1}\vec{g}_{j}(u))_{l}=
\OO(1)
{|\S_j|}{|u|^{-1}} R_{m_l}T_{m_{l},m_{j}}^{-1}
\qquad(l=1,\ldots,j-1)
\Eq(4.3.6)
$$
Moreover, we obtain 
$$
\l\equiv 1-e^{-u}\in\s((1-P_N)^I)
\qquad\Longleftrightarrow\qquad 
G_{j}(u)=0
\Eq(4.3.16)$$
where $G_{j}(u)$ is defined in \eqv(4.3.2e). 
}

\remark
Let us mention that the bound 
on $\Im(u)$ in \eqv(4.3.3) is not optimal and chosen 
just for the sake of convenience. The optimal bounds with 
respect to our control  can easily be derived but 
they are of no particular relevance for the following 
analysis.

\proof
Fix $j=2,\ldots,j_{0}$. Formally we obtain
$$
\KK_{j}(u)^{-1}= 
\left(\1-\DD_(u)^{-1}\NN_j(u)\right)^{-1}\DD_j(u)^{-1}
=
\sum_{s=0}^{\infty}
(\DD_j(u)^{-1}\NN_j(u))^s\DD_j(u)^{-1}
\Eq(4.3.7)
$$
To use these formal calculations and to extract the 
decay estimate 
 in \eqv(4.3.6) 
we must estimate the summands in \eqv(4.3.7). To do this we use a 
straightforward random walk representation for the matrix elements 
$$
\left(\DD_j(u)^{-1}\NN_j(u))^s\DD_j(u)\right)_{l,k}^{-1}=
\sum_{\o:m_l\rightarrow m_{k}\atop |\o|=s}
\prod_{t=1}^{|\o|}
{G_{\o_{t},\S_j}^{\o_{t-1}}(u)
\over
1-G_{\o_{t-1},\S_j}^{\o_{t-1}}(u)}
(1-G_{m_{k},\S_j}^{m_{k}}(u))^{-1},
\quad {1\leq l,k<j}
\Eq(4.3.8)
$$
where  
$\o:m_{l}\rightarrow m_k$ 
denotes a sequence  $\o=(\o_{0},\ldots,\o_{|\o|})$ such 
that $\o_{0}=m_{k}$, $\o_{|\o|}=m_k$, 
$\o_{t}\in \S_{j}\ba (I\cup J)$ 
and $\o_{t-1}\neq\o_{t}$ for all $t=1,\ldots,|\o|$. 
Assuming that the series in \eqv(4.3.7) converges,
\eqv(4.3.8) gives the convenient representation
$$
\eqalign{
(\KK_{j}(u)^{-1}\vec{g}_{j}(u))_{l}=&
\sum_{\o:m_l\rightarrow m_j}
\prod_{t=1}^{|\o|}
{G_{\o_{t},\S_j}^{\o_{t-1}}(u)
\over
1-G_{\o_{t-1},\S_j}^{\o_{t-1}}(u)}
}
\Eq(4.3.8a)
$$
where the sum is now over all walks of arbitrary length. We will 
now show that this sum over random walks does indeed converge 
under our hypothesis.

By virtue of \eqv(4.2.2a) we may apply \eqv(3.7.1) 
for $m\equiv m_l$ and $I\equiv\S_j\ba m_l$ and 
conclude that there are $c>0$ and $C<\infty$ such that 
for all $C^\prime<\infty$ and all $u\in\C$ satisfying 
\eqv(4.3.3)
$$
\eqalign{
G_{m_l,\S_j}^{m_l}(u)-1&=
(1+\OO(\e_N))R_{m_l}^{-1}\left
(u-u_{\S_j\ba m_l}\right)
\left(1+(u-u_{\S_j\ba m_l})\OO(c_N T_{\S_j}) \right)
\cr&
= (1+ \OO(\e_N+2c)) u R_{m_l}^{-1}
}
\Eq(4.3.12)
$$
where we used that $u_{\S_j\ba m_l}\leq c_N \EE_j$. 
In addition, shrinking possibly $c>0$ in \eqv(4.3.3), \eqv(3.4.1) 
implies that for all 
$k,l=1,\ldots,j$, $k\neq l$
$$
G_{m_{k},\S_{j}}^{m_{l}}(u)=
\left(1+\OO(|u|c_N T_{j+1} )\right)
G_{m_{k},\S_{j}}^{m_{l}}(0)\leq \OO(1)
\P[\t^{m_l}_{m_k}\leq \t^{m_l}_{\S_j}]
\Eq(4.3.14)
$$
Using these two bounds, \eqv(4.3.8a) yields
$$
\eqalign{
(\KK_{j}(u)^{-1}\vec{g}_{j}(u))_{l}
\leq  &
\sum_{\o:m_l\rightarrow m_j}
\prod_{t=1}^{|\o|} \OO(1){R_{\o_{t-1}} 
\P[\t^{\o_{t-1}}_{\o_t}\leq \t^{\o_{t-1}}_{\S_j}]}{ |u|^{-1}}
}
\Eq(4.3.10)
$$
To bound the product of probabilities, the following Lemma is useful:

\lemma {\TH(S.3a)}{\it 
Let $\o_0,\o_1,\o_2,\o_k\in\S_j$ such that
$\o_i\neq\o_{i+1}$, for all $i$ and $\o_0\neq\o_k$. Then 
$$
\prod_{t=1}^k  \P[\t^{\o_{t-1}}_{\o_t}\leq \t^{\o_{t-1}}_{\S_j}]  
\leq  \P[
\t^{\o_{0}}_{\o_k}\leq 
\t^{\o_0}_{(\S_j\ba \o_1\ba\dots\ba\o_{k})\cup\o_0}
] 
(\EE_j)^{k-1}
\Eq(4.a.1)
$$
}

\proof The proof is by induction over $k$. For $k=1$ the claim is trivial.
Assume that it for $k=l$. We will show that it holds for $k=l+1$. 
Let $s\equiv\max\{0\leq t\leq l\,|\,\o_t=\o_0\}$. Note that by induction 
hypothesis and definition of $s$,
$$
\prod_{t=s+1}^{l+1}  
\P[\t^{\o_{t-1}}_{\o_t}\leq \t^{\o_{t-1}}_{\S_j}]  
\leq  
\P[\t^{\o_{s}}_{\o_l}\leq 
\t^{\o_s}_{\S_j\ba\o_{s+1}\ba\dots\ba\o_l}]
\P[\t^{\o_l}_{\o_{l+1}}
\leq \t^{\o_l}_{\S_j}] (\EE_j)^{l-s-1}
\Eq(4.a.2)
$$
Now 
$$
\eqalign{
\P[\t^{\o_{s}}_{\o_{l+1}}\leq \t^{\o_s}_{\S_j\ba\o_{s+1}\ba\dots
\ba\o_{l+1}}]&\geq 
\P[\t^{\o_{s}}_{\o_{l+1}}\leq \t^{\o_s}_{\S_j\ba\o_{s+1}\dots\ba\o_{l+1}}
,\t^{\o_s}_{\o_{l}}<
\t^{\o_{s}}_{\o_{l+1}}]\cr
&= \P[\t^{\o_{s}}_{\o_{l}}\leq \t^{\o_s}_{\S_j\ba\o_{s+1}\ba\dots\ba\o_l}]
\P[\t^{\o_l}_{\o_{l+1}}<\t^{\o_l}_{\S_j\ba \o_{s+1}\ba\dots\ba\o_{l+1}}]
\cr
&= \P[\t^{\o_{s}}_{\o_{l}}\leq \t^{\o_s}_{\S_j\ba\o_{s+1}\ba\dots\ba\o_l}]
\frac{\P[\t^{\o_l}_{\o_{l+1}}<
\t^{\o_l}_{(\S_j\ba \o_{s+1}\ba\dots\ba\o_{l-1})\cup\o_{l+1}}]}
{\P[\t^{\o_l}_{(\S_j\ba\o_{s+1}\ba\dots\ba\o_l)\cup\o_{l+1}}<\t^{\o_l}_{\o_l}]}
\cr
&\geq  \P[\t^{\o_{s}}_{\o_{l}}\leq
\t^{\o_s}_{\S_j\ba\o_{s+1}\ba\dots\ba\o_l}]
\frac{\P[\t^{\o_l}_{\o_{l+1}}\leq\t^{\o_l}_{\S_j}]}
{\P[\t^{\o_l}_{(\S_j\ba\o_{s+1}\ba\dots\ba\o_l)\cup\o_{l+1}}<\t^{\o_l}_{\o_l}]}
}
\Eq(4.a.3) 
$$
Now the denominator on the right is,  
$$
\P[\t^{\o_l}_{(\S_j\ba\o_{s+1}\ba\dots\ba\o_l)\cup\o_{l+1}}<\t^{\o_l}_{\o_l}]
\leq \P[\t^{\o_l}_{\S_j\ba\o_l}<\t^{\o_l}_{\o_l}]\leq \EE_j
\Eq(4.a.4)
$$
by \eqv(4.2.2a). Thus, using the obvious bound 
$$
\prod_{t=1}^{s}
\P[\t^{\o_{t-1}}_{\o_t}\leq\t^{\o_{t-1}}_{\S_j}]\leq
(\EE_j)^{s}
\Eq(4.a.5)
$$
and once more that $\o_0\in \S_j\ba\o_{s+1}\ba\ldots\ba\o_{l+1}$, 
\eqv(4.a.3) inserted into \eqv(4.a.2) yields the claim for $k=l+1$
which concludes the proof. \endproof

Using Lemma \thv(S.3a) in \eqv(4.3.14) and the trivial bound $R_{\o_t}\leq 1$,
we get
$$
\eqalign{
(\KK_{j}(u)^{-1}\vec{g}_{j}(u))_{l}&\leq
\P[\t^{m_l}_{m_j}<\t^{m_l}_{m_l}] \sum_{\o:m_l\rightarrow m_j}
\frac{CR_{m_l} }{|u|}
\left(\frac{C\EE_j}{|u|}\right)^{|\o|-1}
\cr&\leq\P[\t^{m_l}_{m_j}<\t^{m_l}_{m_l}] \sum_{k=1}^\infty 
\frac{CR_{m_l}}{|u|}
\left(\frac{C|\S_j|\EE_j}{|u|}\right)^{k-1}
\cr&
\leq \P[\t^{m_l}_{m_j}<\t^{m_l}_{m_l}] \frac {{CR_{m_l} }{|u|^{-1}}}
{1-{C|\S_j|\EE_j}{|u|^{-1}}}
}
\Eq(4.3.15)
$$
If $C|\S_j|\EE |u|^{-1}$ is say smaller than $1/2$, the 
estimate \eqv(4.3.6) follows immediately. 
\eqv(4.3.16) then is a direct consequence of \eqv(3.1.2) and \eqv(4.3.2d), 
since by \eqv(4.3.6) the determinant of $\KK_j(u)$ cannot vanish in
the domain of $u$-values considered.
\endproof

\remark
Defining
$$
\DD_{I}\equiv\op{diag}
(1-G_{m_{l},\MM_{N}}^{m_{l}})_{1\leq l\leq j_{0}},
\qquad
\NN_{I}\equiv \DD_{I}-\GG_{I,\MM_{N}}
\text{and}
(\vec{f}_{I})^t\equiv
(G_{I,\MM_N}^{m_k})_{1\leq k\leq j_{0}}
\Eq(4.r.1)
$$
where $\GG_{I,\MM_{N}}$ is defined in \eqv(3.1.1a), 
a slight modification of the proof above shows that 
for $c>0$ small enough and all 
$\Re(u)<cb_N^{-1}$ such that 
$$
\alpha_{I}\equiv
\min_{m\in \MM_N\ba I}|G_{m,\MM_{N}}^{m}(u)-1|>
(1/c)c_N^{-1}
\max_{m\in\MM_{N}\ba I}T^{-1}_{m,\MM_{N}\ba m}
\Eq(4.r.4)
$$
one can write an absolutely convergent Neumann series for 
$\left(\1-\DD_{I}^{-1}(u)\NN_{I}(u)\right)^{-1}$. 
Furthermore, as a consequence of a random walk 
expansion similar to \eqv(4.3.15) we obtain the bound
$$
(\GG_{I,\MM_{N}}(u)^{-1}\vec{f}_{I}(u))_{l}=
\OO(\alpha_{I}^{-1}c_N^{-1}T_{m_{l},I})
\Eq(4.r.3)
$$
This estimate is needed for the proof of Lemma 5.4.
We are searching for 
solutions $u$ near $u_{\S_{j-1}}$ of the equation 
appearing in \eqv(4.3.16).  
The case $j=1$ is already 
treated in Theorem 3.5. Fix $j=2,\ldots,j_0$. We want to 
apply Lagrange's Theorem to this equation (see [WW]) 
which tells us the following: Fix a point $a\in\C$ and an 
analytic function $\Psi$ defined on a domain containing 
the point $a$. Assume that there is a contour in the 
domain surrounding $a$ such that on this contour the 
estimate $|\Psi(\zeta)|<|\zeta-a|$ holds. Then the equation 
$$
\zeta=a+\Psi(\zeta)
\Eq(4.4.3)
$$
has a unique solution in the interior of the contour. 
Furthermore, the solution can be expanded in the form
$$
\zeta=a+\sum_{n=1}^{\infty}(n!)^{-1}
\partial_{\zeta}^{n-1}\Psi(a)^n
\Eq(4.4.4)
$$
We are in a position to prove

\proposition{\TH(S.4)} {\it
For $j=1,\ldots,j_0$ there is a simple eigenvalue 
$\tilde\l_j=1-e^{-\tilde u_j}<\l_{\S_j}$ such that 
\eqv(4.1.4), \eqv(4.1.3) hold if we replace $\l_j$ by 
$\tilde\l_j$. Let $\tilde\phi_j$ be a corresponding 
eigenfunction. Then \eqv(4.1.6) holds if we replace 
$\phi_{j}$ by $\tilde{\phi}_{j}$.
}

\proof
By means of Theorem 3.5 and \eqv(3.1.3) we may assume that 
$j=2,\ldots,j_0$. The equation in \eqv(4.3.16) 
can be written as
$$
G_{m_{j},\S_{j}}^{m_{j}}(u)-1+\Phi_{j}(\zeta)=0
\Eq(4.4.9)
$$
where we have set 
$\zeta\equiv u\E[\t_{\S_{j-1}}^{m_{j}}]$ and 
$$
\Phi_{j}(\zeta)\equiv
\sum_{l=1}^{j-1}G_{m_{l},\S_{j}}^{m_{j}}(u)
(\KK_{j}(u)^{-1}\vec{g}_{j}(u))_{l}
\Eq(4.4.10)
$$
Fix constants $c>0$, $C<\infty$ and let us 
denote by $U_{j}$ the strip of all $\zeta\in\C$ such that 
$$
T_j/\EE_{j}<\Re(\zeta)<
c T_{j}/T_{j+1},
\qquad
|\Im(\zeta)|< c T_j/(T_{j+1}r_Nc_N)
\Eq(4.4.11)
$$
Putting $\zeta_{\S_{j-1}}\equiv 
u_{\S_{j-1}}\E[\t_{\S_{j-1}}^{m_{j}}]$ it follows 
$\zeta_{\S_{j-1}}=1+ \OO(\e_N)$ from 
\eqv(3.5.0) and \eqv(3.5.4) and we may apply \eqv(3.7.1) 
for $c>0$ small enough and all $\zeta\in U_{j}$ to obtain
$$
G_{m_{j},\S_{j}}^{m_{j}}(u)-1
=
\E[\t_{\S_{j-1}}^{m_{j}}]^{-1}
(1+\OO(\e_N))
R_{m_j}^{-1}
\left(
\zeta-\zeta_{\S_{j-1}}+
(\zeta-\zeta_{\S_{j-1}})^{2}R_{j}(\zeta)
\right)
\Eq(4.4.12)
$$
where 
$\RR_{j}(\zeta)\equiv 
\E[\t_{\S_{j-1}}^{m_{j}}]^{-1}\RR_{\S_{j-1}}(u)$ is defined in 
\eqv(3.7.3). By  \eqv(4.4.12) it follows that 
\eqv(4.4.9) is equivalent to 
$$
\zeta=\zeta_{\S_{j-1}}+\Psi_{j}(\zeta)
\Eq(4.4.15a)
$$
for some function $\Psi_j$ satisfying 
$$
\Psi_{j}(\zeta)=
\E[\t_{\S_{j-1}}^{m_{j}}]
(1+\OO(\e_N))
R_{m_j}^{-1}
\Phi_{j}(\zeta)+
(\zeta-\zeta_{\S_{j-1}})^{2}\RR_{j}(\zeta)
\Eq(4.4.15b)
$$
Using \eqv(L.11) in 
combination with \eqv(4.1.1e), it follows
$$
\RR_{j}(\zeta)=
\OO(T_{j+1}/T_j)
\Eq(4.4.13)
$$
Using \eqv(4.3.6) and the estimate \eqv(4.3.14),  as well as  
\eqv(L.11), we see
that for some $c>0$, $C<\infty$
for all $|\zeta-\zeta_{\S_{j-1}}|\leq1$  
$$
\E[\t_{\S_{j-1}}^{m_{j}}]
\E\left[\t^{m_j}_{m_j}\1_{\t^{m_j}_{m_j}<\t^{m_j}_{\S_j}}\right]
\Phi_{j}(\zeta)=
\sum_{l=1}^{j-1}
\OO\left( c_N^2 T_j^2 T^{-1}_{m_l,m_j} T^{-1}_{m_j,m_l}
\right)\leq \OO(c_N^2 \TT_j^{-1}) 
\Eq(4.4.15)
$$
By means of \eqv(4.4.13) and \eqv(4.4.15) it follows 
for $|\zeta-\zeta_{\S_{j-1}}|\leq1$ 
$$
\Psi_{j}(\zeta)=
\OO(\TT_j^{-1} +T_{j+1}/T_{j})
\Eq(4.4.15c)
$$
Since $\TT_{j}\geq \EE_{j}$, by  \eqv(4.2.2) and 
Definition 1.2, we 
may apply Lagrange's Theorem to \eqv(4.4.15a) giving the 
existence of a solution 
$\tilde\zeta_{j}=\tilde u_{j}\E[\t_{\S_{j-1}}^{m_{j}}]$ 
of \eqv(4.4.9) satisfying 
$|\tilde\zeta_{j}-\zeta_{\S_{j-1}}|<1$.  
We rewrite \eqv(4.4.15a) in the form 
$$
\tilde\zeta_{j}=\zeta_{\S_{j-1}}+
\OO(\TT_j^{-1} +T_{j+1}/T_{j})
\Eq(4.4.16)
$$
By  \eqv(4.3.16) 
$\tilde\l_{j}\equiv 1-e^{\tilde u_{j}}$ defines an 
eigenvalue.
Since from the invertibility of $\KK_j(\tilde{u}_{j})$ it 
follows that the kernel of 
$\GG_j(\tilde u_j)$ is at most one-dimensional, \eqv(3.1.3) 
implies that $\tilde\l_j$ is simple. Using \eqv(3.5.3) and 
\eqv(3.5.4) for $I\equiv\S_{j-1}$, we derive  from \eqv(4.4.16) 
that \eqv(4.1.3) and \eqv(4.1.4) hold, if we replace 
$\l_{j}$ by $\tilde\l_{j}$. Moreover, using 
$\tilde u_{j}<u_{\S_{j}}$ from \eqv(3.1.3), we conclude that 
$$
(\tilde{\phi}_{j}(m_{l}))_{1\leq l<j}=
\tilde{\phi}_{j}(m_{j})
\KK_{j}(\tilde{u}_{j})^{-1}\vec{g}_{j}(\tilde{u}_{j})
\Eq(4.4.16a)
$$
Hence from \eqv(4.3.6) and $\tilde u_{j}=e^{\OO(1)}u_{\S_{j-1}}$ 
we obtain that \eqv(4.1.6) is satisfied if we replace 
$\phi_j$ by $\tilde\phi_j$.
\endproof

Now it is very easy to finish the

\proofof{\thv(S.1)}
Proposition \thv(S.4) tells us that $\l_{k}\leq \tilde{\l}_{k}$ 
for $k=1,\ldots,j_{0}$. Assume now that there is 
$k=2,\ldots,j_0$ such that $\l_{k}<\tilde{\l}_{k}$. 
Let $k=2,\ldots,j_0$ be minimal with this property. Since 
$\tilde{\l}_{k-1}=\l_{k-1}$ is simple, we have 
$\tilde{\l}_{k-1}<\l_k$. 
Lemma \thv(S.3) 
in combination with \eqv(4.3.2e) now 
tells us that for $j=1,\ldots,j_0$ some 
constants $c>0$, $C<\infty$ 
and all $Cc^{-1}_N\EE_j^{-1}<u<c c^{-1}_NT^{-1}_{j+1}$ the 
function $G_{j}(u)$ is strictly monotone decreasing, i.e. 
has at most one zero. 
Hence from \eqv(4.3.16) for $j\equiv k-1$ and 
$G_{k-1}(\tilde{u}_{k-1})=0$ we deduce that 
$u_{k}\geq c c_N^{-1} T^{-1}_k$. But since we already know 
that $u_{k}\leq Cc_N^{-1} T_{k}^{-1}$ for some $C$, 
it then follows from \eqv(4.3.16) for $j\equiv k$ that 
$G_{k}(u_k)=0$ implying the contradiction $\l_k=\tilde{\l}_{k}$. 

Since $\l_{j_0}$ is simple, \eqv(4.3.16) for $j\equiv j_0$ 
and $G_{j_0}(u_{j_0})=0$ 
implies $\l_{j_0+1}>c b_N$, where $c$ denotes 
the constant appearing in \eqv(4.3.3). 

The remaining assertions of Theorem \thv(S.1) then follow 
from Proposition \thv(S.4). 
\endproof

\medskip

\chap{6. The distribution function}6

The objective of this chapter is to show how the 
structure of the low lying spectrum implies a 
precise control of the  
distribution 
function of the times $\t^{m}_I$, in cases where
Theorem \thv(LL.5) applies, i.e. 
$I\sb\MM_{N}$, $I,\MM_{N}\ba I\neq\em$, and 
$m_1\in\MM_{N}\ba I$,  $T_{I}=T_{m_1,I}$. 
It is already known that the normalized distribution 
function converges weakly to the exponential 
distribution (see [BEGK] for the sharpest estimates 
beyond weak convergence in the most general case). 

The proof of these results proceeds by inverting the 
 Laplace transforms $G^m_I(u)$, making use of the information about the
analytic structure of these functions that is contained in the
spectral decomposition of the 
low lying spectrum of $(1-P_{N})^{I}$ obtained in the previous section.

Let us denote by $\LL_N$ the Laplace transform of the 
complementary distribution function, i.e.
$$
\LL_N(u)\equiv\LL_{N,I}^{m_1}(u)\equiv
\sum_{t=0}^\infty e^{ut}\P[\t_I^{m_1}>t]
\qquad(\Re(u)<u_{I}),
\Eq(5.2.1)$$
where $u_{I}$ is defined in \eqv(3.5.0). The 
Perron-Frobenius Theorem gives 
$\lim(1/t)\log\P[\t_I^{m_1}>t]=-u_I$. Hence the Laplace 
transform defined above is
locally uniformly exponentially convergent. In order 
to obtain the continuation 
of $\LL_{N}$ to the whole plane we perform a 
partial summation in the sum on the right-hand side 
of \eqv(5.2.1) and get
$$
\LL_{N}(u)=
{G_{I,I}^{m_1}(u)-1 \over e^u-1}.
\Eq(5.2.2)$$
Invoking \eqv(1.5.1) a straightforward computation for 
$\l\equiv 1-e^{-u}$ shows that
$$
G_{I,I}^{x}(u)=
((1-P_{N})^{I}-\l)^{-1}(\1_{I^{c}}P_{N}\1_{I})(x)
\qquad(x\notin I),
\Eq(5.2.3)$$ 
Hence $\LL_{N}$ is a meromorphic function with poles in 
$u\in\{u_{1},\ldots,u_{|\G_{N}\ba I|}\}$, where 
we recall the definition of the eigenvalues 
$\l_j=1-e^{-u_j}$ for $j=1,\ldots,|\G_{N}\ba I|$ 
prior to Theorem 4.1. Since $\LL_{N}$ 
is $2\pi$-periodic in the imaginary direction, a short 
computation yields
$$
\P[\t_I^{m_1}>t]={1\over 2\pi i}
\int_{-i\pi}^{i\pi}e^{-tu}\LL_{N}(u)du.
\Eq(5.1.18)
$$
Deforming the contour in \eqv(5.1.18) gives 
for $u_{j_0}<\a<u_{j_{0}+1}$ and 
$U_{\a}\equiv(0,\a)\times(-\pi,\pi)$
$$
\P[\t_I^{m_1}>t]=
{1\over 2\pi i}\int_{\a-i\pi}^{\a+i\pi}
e^{-tu}\LL_{N}(u)du
-\sum_{u_{j}\in U_{\alpha}}e^{-tu_{j}}\op{res}_{u_{j}}
\LL_{N},
\Eq(5.1.19)
$$
where $\op{res}_{u}\LL_N$ denotes the residue of $\LL_N$ 
at $u$. Here we have used that periodicity of $\LL_{N}$ 
shows that the integrals over $[\a+i\pi,i\pi]$ and 
$[-i\pi,\a-i\pi]$ cancel and that the poles $u_{j}$, 
$j=1,\ldots,j_{0}$, are simple.

Our main result can be formulated as follows:

\theo{\TH(LA.1)} {\it
Let $j_0\equiv |\MM_N\ba I|$. There is $c>0$ such that 
for some $c>0$, 
$$
\eqalign{
\P[\t_I^{m_1}>t] = 
-\sum_{j=1}^{j_0}
e^{-tu_j}\op{res}_{u_j}\LL_N+
e^{-tc b_N^{-1}}(2\pi i)^{-1}
\int_{-i\pi}^{i\pi}
e^{-tu}\LL_N(u)du,
}
\Eq(5.1.2)
$$
where  $u_j=-\ln (\l_I+1)$ and $\l_j$ are  the eigenvalues
of $(1-P_N)^I$ that are estimated in Theorem \thv(S.1). 
Moreover, the residues satisfy 
$$
\op{res}_{u_{1}}\LL_{N}=
-1+\OO\left(R_{m_1}c_N T_2/T_1\right),
\qquad
\op{res}_{u_{j}}\LL_{N}=
\OO\left(R_{m_1}c_N T_{j}/T_{1}\right)
\qquad(j=2,\ldots,j_0)
\Eq(5.1.3)$$
while the remainder integral on the right-hand side 
of \eqv(5.1.2) is bounded by 
$$
(2\pi i)^{-1}\int_{-i\pi}^{i\pi}
e^{-tu}\LL_N(u)du=
\OO\left(c_N^{-1}b_N^{-2}|\G_N|^2/T_1\right).
\Eq(5.1.3b)
$$
}

\remark 
Recalling  \eqv(L.11) and Theorem 4.1,
one sees that Theorem \thv(LA.1) implies that the 
distribution of $t^{m_1}_I$ is to a remarkable precision a 
pure exponential.
In particular, one has the

\corollary{\TH(LA.2)} {\it
Uniformly in $t\in \E[\t_I^{m_1}]^{-1}\N$
$$
\P[\t_I^{m_1}>t\E[\t_I^{m_1}]]
=
\left(1+\OO(R_{m_1}c_N T_2/T_1)
\right)
e^{-t\left(
1+\OO(R_{m_1}c_N T_2/T_1)\right)}.
\Eq(5.1.3a)
$$
}

We start with the computation of the residue of the 
Laplace transform at $u_1$.

\lemma{\TH(LA.3)} {\it
$$
\op{res}_{u_{1}}\LL_{N}=
-1+\OO(R_{m_1}c_N T_{2}/T_{1}).
\Eq(5.2.4)
$$
}

\proof
From \eqv(3.5.1) for $m\equiv m_1$ and the renewal 
relation \eqv(1.6.1) and \eqv(5.2.2)  follows
$$
\op{res}_{u_{1}}\LL_{N}=
\lim_{u\rightarrow u_{1}}
{G_{I,m_1}^{m_{1}}(u)\over e^{u}-1}
{u-u_{1}\over 
G_{m_1,I}^{m_1}(u_{1})-G_{m_1,I}^{m_1}(u)}
= -{1\over e^{u_{1}}-1}
{G_{I,m_1}^{m_{1}}(u_{1})\over
\del_u {G}_{m_1,I}^{m_1}(u_{1})}.
\Eq(5.2.5)$$
Since $u_1=e^{\OO(1)}N^{-1}R_{m_1}T^{-1}_1$, \eqv(3.4.1) for $k=0,1$ 
gives for some $C<\infty$
$$
{G_{I,m_1}^{m_1}(u_{1})\over
\del_u {G}_{m_1,I}^{m_1}(u_{1})}
=(1+\OO( R_{m_1} c_N T_{2}/T_{2}))
{G_{I,m_1}^{m_1}(0)\over\del_u {G}_{m_1,I}^{m_1}(0)}.
\Eq(5.2.6)$$
Hence \eqv(5.2.4) follows from 
\eqv(5.2.5) in combination with \eqv(4.1.3) and 
\eqv(2.6.1). 
\endproof

In general we cannot prove lower bounds for 
the higher residues for the reason described in the 
remark after Theorem 4.1. However, we can show that they are 
very small: 

\lemma{\TH(LA.4)} {\it
$$
\op{res}_{u_{j}}\LL_{N}=
\OO\left(T_{j}/T_{1})\right)
\qquad(j=2,\ldots,j_0).
\Eq(5.3.1)
$$
}

\proof
For fixed $j=0,\ldots,j_0$ we compute, using 
\eqv(5.2.2) and \eqv(5.2.3),
$$
\eqalign{
\op{res}_{u_{j}}\LL_{N}
= &
\lim_{u\rightarrow u_{j}}
{1\over e^u-1}
{u-u_{j}\over (1-e^{-u_{j}})-(1-e^{-u})}
{\la \1_{I^{c}}P_{N}\1_{I},\phi_{j}\ra_{\Q_{N}}
\over
(||\phi_{j}||_{\Q_{N}})^2}\phi_{j}(m_1)
\cr = & 
-{e^{u_{j}}\over e^{u_{j}}-1}
{\la \1_{I^{c}}P_{N}\1_{I},\phi_{j}\ra_{\Q_{N}}
\over
(||\phi_{j}||_{\Q_{N}})^2}\phi_{j}(m_1).
}
\Eq(5.3.2)
$$
We can assume that $\phi_j(m_j)=1$.
We can express $\phi_j(x)$, using the definition  
\eqv(3.1.3),  Lemma \thv(C.4), and Theorem \thv(S.1)  
in the form 
$$
\eqalign{
\phi_j(x)=&(1+\OO(\g))K_{m_j,\S_j}^{x}(0)
+\sum_{l=1}^{j-1}\OO( T_j/T_{m_l,m_j})
(1+\OO(\g))K_{m_l,\S_j}^{x}(0)
\cr = &
(1+\OO(\g))\P[\s_{m_j}^{x}<\t_{\S_{j-1}}^{x}]+
\OO(\g).
}
\Eq(5.3.11a)
$$
where  $\g\equiv R_{m_j}\max (\TT^{-1},T_{j+1}/T_{j})$.
Using Lemma \thv(LL.3), one sees easily that this implies that for any $\e>0$,
$$
(||\phi_j||_{\Q_N})^2
\geq (1+\OO(e^{-N\g}))
\Q_N(\{x\in\G_N\,|\,|x-m_j|<\vep/2\})
\geq (1-\e) \QQ_N(A(m_j))
\Eq(5.3.11)
$$
From \eqv(3.1.3) we conclude that,  for $J\equiv\S_{j}$,
$$
\eqalign{
\la \1_{I^{c}}P_{N}\1_{I},\phi_{j}\ra_{\Q_{N}}
=&
\sum_{k=1}^{j}\phi_{j}(m_k)
\sum_{x\in\G_{N}\atop y\in I}
\Q_{N}(x)P_{N}(x,y)
K_{m_k,\S_j}^{x}(u_{j})
\cr = &
\sum_{k=1}^{j}\phi_{j}(m_k)
\sum_{x\in\G_{N}\atop y\in I}
\Q_{N}(y)P_{N}(y,x)K_{m_k,\S_j}^{x}(u_{j}),
}
\Eq(5.3.2a)
$$
where we have used the symmetry of $P_N$. 
Applying \eqv(1.5.1) and 
\eqv(1.7.1) to the right-hand side of \eqv(5.3.2a) 
we get
$$
\eqalign{
\la \1_{I^{c}}P_{N}\1_{I},\phi_{j}\ra_{\Q_{N}} = &
\sum_{k=1}^{j}\phi_{j}(m_k)
\sum_{y\in I}\Q_{N}(y)G_{m_k,\S_j}^{y}(u_{j})
\cr = &
\sum_{k=1}^{j}\phi_{j}(m_k)
\Q_{N}(m_k)G_{I,\S_j}^{m_k}(u_{j}).
}
\Eq(5.3.3)
$$
Using that $\phi_j(m_j)=1$, we deduce from \eqv(4.1.6) and reversibility that 
$$
\Q_{N}(m_k)\phi_j(m_k)=\Q_N(m_j)
\OO( R^{-1}_{m_j} T_j/T_{m_j,m_k})
\Eq(5.3.3a)
$$
Combining   \eqv(5.3.3a) with \eqv(4.3.14), 
\eqv(5.3.11), and, once more,  \eqv(4.1.6) with $k\equiv 1$, 
gives
$$
\eqalign{
(||\phi_j||_{\Q_N})^{-2}\phi_{j}(m_1)
\la \1_{I^{c}}P_{N}\1_{I},\phi_{j}\ra_{\Q_{N}}
=&
\sum_{k=1}^{j}
\OO\left(R_{m_j}
\frac{T_j^2}{ T_{m_1,m_j}T_{m_j,m_k}T_{m_k,I}}
\right)
\cr = &
\OO\left(R_{m_j}
\frac{T_j^2}{T_{m_1,m_j}T_{m_j,I}}
\right),
}
\Eq(5.3.8)$$
where we have used Lemma \thv(S.3a) for the sequences 
$\o=(m_j,m_k,m)$ in the 
last equation. 
It is easy to verify that 
$$
\frac{T_j^2}{T_{m_1,m_j}T_{m_j,I}}
 \leq  
\frac{T_j}{T_{m_j,I\cup m_1}T_1}.
\Eq(5.3.10a)
$$
Inserting \eqv(5.3.8) and \eqv(5.3.10a) into 
\eqv(5.3.2), using $u_j=R_{m_j} T_j^{-1}(1+o(1))$ 
 and $T_{m_j,I\cup m_1}\geq T_j$, we arrive at 
\eqv(5.3.1).
\endproof

The last ingredient for the proof of Theorem 
\thv(LA.1) consists in estimating of the remainder 
integral in \eqv(5.1.2). This essentially boils  
down to 

\lemma{\TH(LA.5)} {\it
There is $\d>0$ such that for all 
$\d^{-1}R_{m_1}T_{j_{0}}<\alpha<\d b_N|\G_N|^{-1}$ 
and 
all $\l\equiv 1-e^{-u}$ on the circle 
$|\l-1|=e^{-\alpha}$ we have
$$
G_{I,I}^{m_1}(u)=\OO(\alpha^{-1}c_N^{-1}T_1^{-1}).
\Eq(5.4.1)$$
}

\proof
From the strong Markov property \eqv(1.3.6) for 
$J\equiv I$ and $L\equiv\MM_N\ba I$ we obtain for 
$\Re(u)<u_{\MM_N}$ 
$$
K_{I,I}^{x}(u)=K_{I,\MM_N}^x(u)+
\sum_{l=1}^{j_0} 
K_{I,I}^{m_l}(u)K_{m_l,\MM_N}^x(u)
\qquad(x\in\G_{N}).
\Eq(5.1.6)$$
Applying $(1-P_N-\l)^I$ to both sides of the previous 
equation and evaluating the resulting equation at 
$x=m_k$, $k=1,\ldots,j_0$, we conclude, as in 
\eqv(3.1.6), via \eqv(1.5.2) and \eqv(1.3.4) that 
$$
0=
-G_{I,\MM_N}^{m_k}(u)+
\sum_{l=1}^{j_0}
G_{I,I}^{m_l}(u)(\d_{lk}-G_{m_l,\MM_N}^{m_k}(u)).
\Eq(5.1.7)$$
Thus the vector
$$
\vec{\psi}_{\l}\equiv
(G_{I,I}^{m_l}(u))_{1\leq l\leq j_{0}}
\Eq(5.1.8)$$
solves the system of equations
$$
\GG_{I,\MM_{N}}(u)\vec{\psi}_{\l}=
\vec{f}_{I}(u),
\Eq(5.1.9)$$
where $\GG_{I,\MM_{N}}(u)$ and $\vec{f}_{I}(u)$ are defined 
in \eqv(3.1.1a) and \eqv(4.r.1), respectively. 
In order to be able to apply \eqv(4.r.3) we claim that 
for some $\d,c>0$, for all $u=\alpha+iv$, 
$v\in [-\pi,\pi]$, and for all $m\in\MM_N\ba I$
$$
|G_{m,\MM_N}^{m}(u)-1|\geq c\alpha.
\Eq(5.1.10)$$
We first observe that \eqv(1.3.1) shows that, for all 
$\Re(u^\prime)<u_{\MM_{N}}$,
$$
\Q_N(m)(G_{m,\MM_N}^{m}(u)-1)=
-e^{u}\la ((1-P_N)^{\MM_N\ba m}-\l)
K_{m,\MM_N}^{(\cdot)}(u),
K_{m,\MM_N}^{(\cdot)}(u^\prime) \ra_{\Q_N},
\Eq(5.1.11)$$
where we have extended the inner product to $\C^{\G_N}$ 
in the canonical way such that it is $\C$-linear in the 
second argument. For $|v\pm\pi|\leq \pi/3$ we simply get 
from \eqv(5.1.11), for $u^\prime\equiv u$ and some $c>0$, 
using that $\s((1-P_{N})^{\MM_{N}\ba m})\sb(0,1)$, 
$$
\eqalign{
|\Q_N(m) & \op{Re}(e^{-u}(G_{m,\MM_N}^{m}(u)-1))|
\cr = & \left|\left\la
((1-P_{N})^{\MM_{N}\ba m}-(1+e^{-\alpha}|\cos(v)|)
K_{m,\MM_N}^{(\cdot)}(u),K_{m,\MM_N}^{(\cdot)}(u)
\right\ra_{\Q_{N}}\right|
\cr \geq &
(1+ce^{-\alpha}-1)
(||K_{m,\MM_N}^{(\cdot)}(u)||_{\Q_N})^2
\cr \geq &
c e^{-\alpha}\Q_N(m).
}
\Eq(5.1.12)$$
For $|v+\pi|>\pi/3$, $|v-\pi|>\pi/3$ and 
$|v|>\alpha$, we derive from \eqv(5.1.11) for 
$u^\prime\equiv u$ and some $c>0$ 
$$
\eqalign{
|\Q_N(m)\op{Im}(e^{-u}(G_{m,\MM_N}^{m}(u)-1))|
= &
|\sin(v)|e^{-\alpha}
(||K_{m,\MM_N}^{(\cdot)}(u)||_{\Q_N})^2
\cr \geq &
\Q_N(m)c \alpha e^{-\alpha}.
}
\Eq(5.1.13)$$
In the remaining case, namely where $|v|\leq\alpha$, 
we use \eqv(5.1.11) for $u^\prime\equiv u_{\MM_{N}\ba m}$ 
and obtain via \eqv(3.1.3), for $I\equiv\MM_N\ba m$, $J\equiv m$, that
$$
|\Q_N(m)e^{-u}(G_{m,\MM_N}^{m}(u)-1)|=
|\bar{\l}-\l_{\MM_{N}\ba m}|\,
|\la K_{m,\MM_{N}}^{(\cdot)}(u),
K_{m,\MM_{N}}^{(\cdot)}(u_{\MM_{N}\ba m})\ra_{\Q_{N}}|.
\Eq(5.1.14)$$
From \eqv(3.4.1) it follows for some $\d>0$ uniformly in 
$x\in\G_{N}$ and $|v|\leq \alpha$
$$
K_{m,\MM_{N}}^{x}(u)=
\left(1+\d\OO(1)\right)
K_{m,\MM_{N}}^{x}(u_{\MM_{N}\ba m}).
\Eq(5.1.15)$$
Since the minimum of the function 
$|\bar{\l}-\l_{\MM_N\ba m}|$ is attained at 
$\l=1-e^{-\alpha}$, we conclude from \eqv(5.1.14) and 
\eqv(5.1.15) in combination with \eqv(3.1.3) for 
$J\equiv m_{1}$ and 
\eqv(5.3.11) for some $c>0$ and all 
$|v|\leq \alpha$ that  
$$
\eqalign{
|\Q_N(m)e^{-u}(G_{m,\MM_N}^{m}(u)-1)|\geq &
c|\bar{\l}-\l_{\MM_N\ba m}|
(||K_{m,\MM_N}^{(\cdot)}(u_{\MM_{N}\ba m})||_{\Q_N})^{2}
\cr \geq & c^2\Q_N(A(m))(1-e^{-\alpha}).
}
\Eq(5.1.16)
$$
\eqv(5.1.16), \eqv(5.1.13) and \eqv(5.1.12) prove 
\eqv(5.1.10). Since by definition \eqv(4.1.1d) 
and \eqv(4.2.2) it follows that 
$$
T_{j_{0}}=T_{m_{j_0},\MM_N\ba m_{j_0}}=
\min_{m\in\MM_N}T_{m,\MM_N\ba m}\geq b^{-1}_N,
\Eq(5.1.20)
$$
$b_N$ is defined 
in Definition 1.1, 
combining \eqv(5.1.10) with 
\eqv(4.r.3) shows that the solution of \eqv(5.1.9) 
satisfies  
$$
\psi_{\l}(m_1) = (\vec{\psi}_{\l})_1 =
\OO\left(\alpha^{-1} c_N^{-1}/T_1\right).
\Eq(5.1.17)
$$

\proofof{ Theorem \thv(LA.1)} The proof of Theorem \thv(LA.1) now is 
reduced to 
the  application of the Laplace inversion formula 
and estimation of the remainder integral.
In view of \eqv(5.2.4) and \eqv(5.3.1) it remains to 
estimate the remainder integral on the right-hand side 
of \eqv(5.1.19). But this is by means of \eqv(5.2.2) 
and \eqv(5.2.3) in combination with \eqv(5.4.1) for 
$\alpha\equiv cb_N|\G_N|^{-1}$, $0<c<\d$, fairly easy.
\endproof

\bigskip

\chap{References}0

\item{[BBG]} G. Ben Arous, A. Bovier, and V. Gayrard, ``Aging in the random 
energy model under Glauber dynamics'', in preparation (2000).
\item{[BEGK]} A. Bovier, M.Eckhoff, V. Gayrard, and M. Klein, ``Metastability in Stochastic 
Dynamics of Disordered Mean-Field Models``, to appear in 
Probab. Theor. Rel. Fields (2000).
\item{[BK]} G. Biroli and J. Kurchan, ``Metastable states in glassy systems'',
\hfill\break{\tt http://www.xxx.lanl.gov/cond-mat/0005499} (2000).
\item{[DS]} P.G. Doyle and J.L. Snell, ``Random walks and electrical
networks'',  
Carus Mathematical Monographs, 22,
Mathematical Association of America, Washington, DC, 1984.
\item{[DV]} M.D.  Donsker and S.R.S. Varadhan,  ``On the principal eigenvalue 
of 
second-order elliptic differential
operators'', 
Comm. Pure Appl. Math. {\bf 29},  595-621 (1976).  
\item{[EK]} M. Eckhoff and M. Klein, ``Metastability and low lying spectra
in non-reversible Markov chains'', in preparation (2000).
\item{[FW]} M.I. Freidlin and A.D. Wentzell, ``Random perturbations of 
dynamical systems'', Springer, Berlin-Heidelberg-New York, 1984.
\item{[GM]} B. Gaveau and M. Moreau, ``Metastable relaxation times and absorbtion
probabilities for multidimensional stochastic systems'', J. Phys. A: Math. Gen.
{\bf 33}, 4837-4850 (2000).
\item{[GS]}  B. Gaveau and L.S. Schulman, ``Theory of nonequilibrium 
first-order phase transitions for stochastic dynamics'',  J.
Math. Phys. {\bf 39}, 1517-1533 (1998
\item{[Li]} T.M. Liggett, ``Interacting particle systems'', Springer, Berlin, 
1985.
\item{[M]} P. Mathieu, ``Spectra, exit times and long times asymptotics in 
the zero white noise limit'', Stoch. Stoch. Rep. {\bf 55}, 1-20 (1995).
\item{[S]} Ch. Sch\"utte, ``Conformational dynamics: modelling, theory, 
algorithm, and application to biomolecules'', preprint SC 99-18, 
ZIB-Berlin (1999).
\item{[SFHD]} Ch. Sch\"utte, A. Fischer, W.  Huisinga, and P. Deuflhard, ``A 
direct approach to conformational dynamics based on
hybrid Monte Carlo'', J. Comput. Phys. {\bf 151}, 146-168 (1999).
\item{[Sc]} E. Scoppola, ``Renormalization and graph methods for Markov 
chains'', Advances in dynamical systems and quantum
physics (Capri, 1993), 260-281, World Sci. Publishing, River Edge, NJ, 1995.
\item{[W]}  A.D. Wentzell, ``On the asymptotic behaviour of the greatest
 eigenvalue of a second order elliptic differential operator with a small 
parameter  in the higher derivatives'', Soviet Math. Docl. {\bf 13}, 13-17 
(1972).
\item{[WW]} E.T. Whittaker and G.N. Watson, ''A course of Modern Analysis'', 
Cambridge University Press, (1958).

\end